\documentclass[english,a4paper,10pt]{article}

\usepackage{xcolor}
\usepackage{amsfonts,amsmath,amstext,amsbsy,amsthm,amscd}
\usepackage{graphics}
\usepackage{epsfig}
\usepackage{geometry}
\usepackage{amssymb}
\usepackage[latin1]{inputenc}
\usepackage{babel,indentfirst,amsfonts,array}

\newtheorem{thm}{Theorem}[section]
\newtheorem{lemme}[thm]{Lemma}
\newtheorem{prop}[thm]{Proposition}
\newtheorem{cor}[thm]{Corollary}

\newtheorem{defi}[thm]{Definition}

\newtheorem{hyp}[thm]{Hypothesis}

\def\T{\mathbb T}
\def\R{\mathbb R}

\def\Q{\mathbb Q}
\def\Z{\mathbb Z}
\def\N{\mathbb N}
\def\({\left(}
\def\){\right)}
\def\[{\left[}
\def\]{\right]}
\def\fin{\hfill\square}

\def\fin{\hfill $\square$}

\title{\textsc{Planar random walk in a stratified quasi-periodic environment}
\author{Julien Br\'emont}
\date{Universit\'e Paris-Est Cr\'eteil,~d\'ecembre 2020}
}

\begin{document}

\maketitle

\setcounter{page}{1}

\begin{abstract}
Completing former works \cite{jb1,jb2,jb3}, we study the recurrence of inhomogeneous Markov chains in the plane, when the environment is horizontally stratified and the heterogeneity of quasi-periodic type.

\end{abstract}

\footnote{
\begin{tabular}{l}\textit{AMS $2020$ subject classifications~: 37E05, 37E10, 60G17, 60J10, 60K37.} \\
\textit{Key words and phrases~: Markov chain, recurrence, stratified environment, quasi-periodic environment.} 
\end{tabular}}

\section{Introduction}
This article investigates the question of the recurrence of a class of inhomogeneous Markov chains in the plane, assuming the environment invariant under horizontal translations. This type of random walks were first considered by Matheron and de Marsily \cite{mama} around 1980, with a motivation coming from hydrology and the modelization of pollutants diffusion in a porous and stratified ground. In 2003, a discrete version was introduced by Campanino and Petritis in \cite{cp}.

\medskip
As in \cite{jb1,jb2,jb3}, we consider an extension of the latter, restricting here to the plane and simplifying a little the hypotheses. We shall define a Markov chain $(S_k)_{k\geq0}$ in $\Z^{2}$, starting at the origin, such that the transition laws are constant on each stratum $\Z\times\{n\}$, $n\in \Z$. The first and second coordinates will be respectively called ``horizontal" and ``vertical". For each (vertical) $n\in\Z$, let positive reals $\alpha_n,\beta_n,\gamma_n$, with $\alpha_n+\beta_n+\gamma_n=1$, and a probability measure $\mu_n$ so that:

\begin{hyp}
\label{hypo} 

$ $

\noindent
$\exists\eta>0$, $\forall n\in\Z$, $\min\{\alpha_n,\beta_n,\gamma_n\}\geq\eta$, $\mbox{Supp}(\mu_n)\subset \Z\cap]-1/\eta,1/\eta[$, $\mu_n(0)\leq 1-\eta$.
\end{hyp}

\medskip
\noindent
The transition laws of $(S_k)_{k\geq0}$ are defined, for all $(m,n)\in\Z^2$ and $ r\in\Z$, by:

$$(m,n)\overset{\alpha_n}{\longrightarrow}(m,n+1),~~~(m,n)\overset{\beta_n}{\longrightarrow}(m,n-1),~~(m,n)\overset{\gamma_n\mu_n(r)}{\longrightarrow}(m+r,n).$$

\medskip
\noindent
Here is the corresponding picture:

\begin{center}
\resizebox{0.4\linewidth}{!}{\begin{picture}(0,0)%
\includegraphics{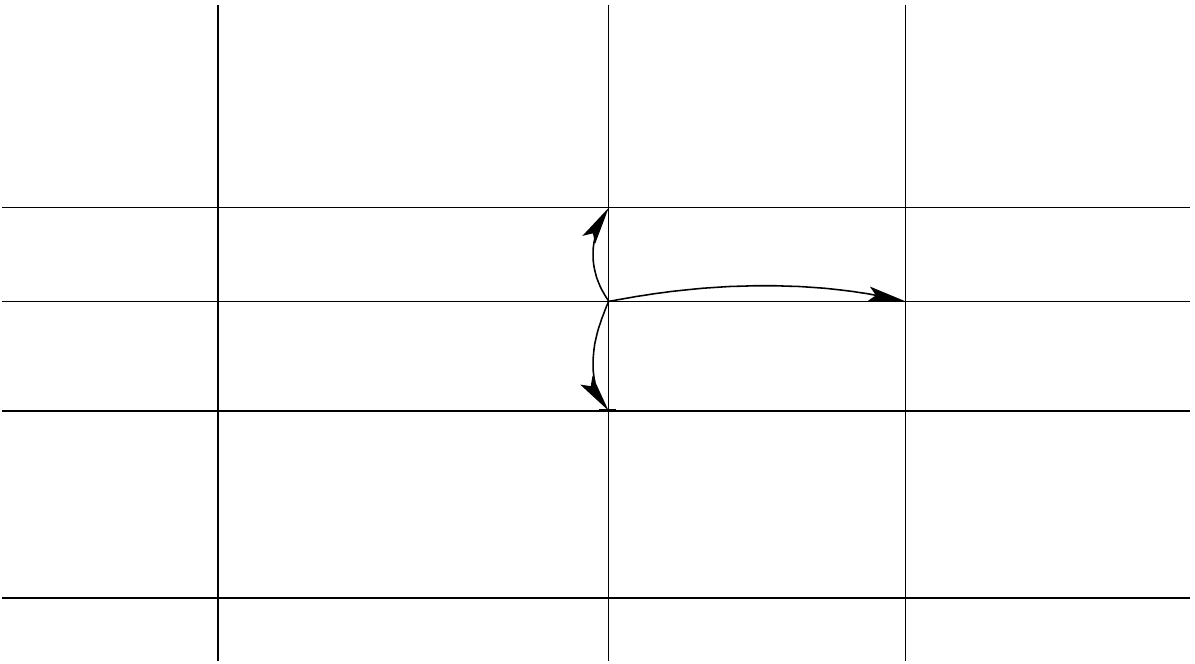}%
\end{picture}%
\setlength{\unitlength}{3947sp}%
\begingroup\makeatletter\ifx\SetFigFont\undefined%
\gdef\SetFigFont#1#2#3#4#5{%
  \reset@font\fontsize{#1}{#2pt}%
  \fontfamily{#3}\fontseries{#4}\fontshape{#5}%
  \selectfont}%
\fi\endgroup%
\begin{picture}(5724,3174)(4279,-4873)
\put(8026,-4786){\makebox(0,0)[lb]{\smash{{\SetFigFont{12}{14.4}{\rmdefault}{\mddefault}{\updefault}{\color[rgb]{0,0,0}$m+r$}%
}}}}
\put(4951,-3361){\makebox(0,0)[lb]{\smash{{\SetFigFont{12}{14.4}{\rmdefault}{\mddefault}{\updefault}{\color[rgb]{0,0,0}$n$}%
}}}}
\put(4801,-2611){\makebox(0,0)[lb]{\smash{{\SetFigFont{12}{14.4}{\rmdefault}{\mddefault}{\updefault}{\color[rgb]{0,0,0}$n+1$}%
}}}}
\put(4801,-3886){\makebox(0,0)[lb]{\smash{{\SetFigFont{12}{14.4}{\rmdefault}{\mddefault}{\updefault}{\color[rgb]{0,0,0}$n-1$}%
}}}}
\put(5476,-1861){\makebox(0,0)[lb]{\smash{{\SetFigFont{12}{14.4}{\rmdefault}{\mddefault}{\updefault}{\color[rgb]{0,0,0}$\Z$}%
}}}}
\put(6826,-4786){\makebox(0,0)[lb]{\smash{{\SetFigFont{12}{14.4}{\rmdefault}{\mddefault}{\updefault}{\color[rgb]{0,0,0}$m$}%
}}}}
\put(9676,-4486){\makebox(0,0)[lb]{\smash{{\SetFigFont{12}{14.4}{\rmdefault}{\mddefault}{\updefault}{\color[rgb]{0,0,0}$\Z$}%
}}}}
\put(5026,-4786){\makebox(0,0)[lb]{\smash{{\SetFigFont{12}{14.4}{\rmdefault}{\mddefault}{\updefault}{\color[rgb]{0,0,0}$0$}%
}}}}
\put(6751,-3436){\makebox(0,0)[lb]{\smash{{\SetFigFont{12}{14.4}{\rmdefault}{\mddefault}{\updefault}{\color[rgb]{0,0,0}$\beta_n$}%
}}}}
\put(6901,-2911){\makebox(0,0)[b]{\smash{{\SetFigFont{12}{14.4}{\rmdefault}{\mddefault}{\updefault}{\color[rgb]{0,0,0}          $\alpha_n$}%
}}}}
\put(7501,-2911){\makebox(0,0)[lb]{\smash{{\SetFigFont{12}{14.4}{\rmdefault}{\mddefault}{\updefault}{\color[rgb]{0,0,0}$\gamma_n\mu_n(r)$}%
}}}}
\end{picture}%
}
\end{center}

\medskip
\noindent
The family of transition laws is called the environment, here identified to $((\alpha_n,\beta_n,\gamma_n,\mu_n))_{n\in\Z}$. Introduce the local horizontal drift $\varepsilon_n:=\sum_{r\in\Z}r\mu_n(r)$ at each vertical $n$, i.e. the expectation of $\mu_n$. The special case when $p_n=q_n$, $n\in\Z$, is called the ``vertically flat model". 

\medskip
With respect to a fully inhomogeneous random walk in the plane, the horizontal stratification of the environment brings the notable simplification that the vertical component of $(S_k)_{k\geq0}$, in restriction to vertical jumps, is a one-dimensional Markov chain. We call it the ``vertical random walk''. With $\alpha'_n=\alpha_n/(\alpha_n+\beta_n)$ and $\beta'_n=\beta_n/(\alpha_n+\beta_n)$, its transition laws on $\Z$ are:

$$n\overset{\alpha'_n}{\longrightarrow}n+1,~~~n\overset{\beta'_n}{\longrightarrow}n-1.$$ 

\smallskip
\noindent
From this, the model inherits some ``product structure''. For instance, for $(S_k)_{k\geq0}$ to be recurrent itself, the vertical random walk has first to be. The conditions for this are known for a long time in the context of birth and death processes (cf Karlin and McGregor \cite{KMG}). Placing in this case, the study of the recurrence of $(S_n)$ essentially reduces to the analysis of the horizontal displacement.  

\medskip

Let us discuss former results concerning the recurrence/transience of such a Markov chain, first for the vertically flat model. The vertical random walk is then simple random walk on $\Z$, which is well-known to be recurrent. The main object of study has been the Campanino-Petritis model \cite{cp}, corresponding to $\alpha_n=\beta_n=\gamma_n=1/3$ and $\mu_n=\delta_{x_n}$, for some $(x_n)_{n\in\Z}\in\{\pm 1\}^{\Z}$, where $\delta_x$ is Dirac measure at $x$. Recurrence is shown in \cite{cp} when $x_n=(-1)^n$ and transience when $x_n=1_{n\geq 1}-1_{n\leq0}$ or if $(x_n)$ are typical realizations of independent and identically distributed ($i.i.d.$) random variables with law $(\delta_1+\delta_{-1})/2$. For random $(x_n)$, transience results were obtained by Guillotin-Plantard and Le Ny \cite{GPLN},
when the $(x_n)$ are independent, but with different marginals, and also by P\`ene \cite{Pene}, under some hypotheses of stationarity and decorrelation. Devulder and P\`ene \cite{DP} showed transience for a more general vertically flat model, where $\alpha_n=\beta_n$, the $(\gamma_n)$ are $i.i.d.$ non-constant and $\mu_n=\delta_{x_n}$, for an arbitrary $(x_n)\in\{\pm 1\}^{\Z}$. In \cite{cp2}, for their initial model, Campanino and Petritis studied the case of a random perturbation of a periodic $(x_n)$. 

\medskip
\noindent
In \cite{jb1}, for the general vertically flat case, a complete recurrence criterion was given. The asymptotic behaviour of the random walk is governed by the sums $(\gamma_{-m}\varepsilon_{-m}/\alpha_{-m}+\cdots+\gamma_{n-1}\varepsilon_{n-1}/\alpha_{n-1})_{-m\leq 0\leq n}$, associated with some horizontal flow defined by the environment and transverse to the vertical layer $[-m,n)$. The central role is played by a two-variable function $\Phi(-m,n)$, introduced below, measuring the ``horizontal dispersion'' of the previous flow between vertical levels $-m$ and $n$. The quantity deciding for the recurrence/transience of $(S_k)_{k\geq0}$ computes some ``capacity of dispersion to infinity'' of the environment. The abstract form of the criterion in \cite{jb1} is directly extracted from a Poisson kernel in a half-plane. It seems to be related to some notion of curvature at infinity of the level lines of the function $\Phi(-m,n)$. Several examples were next presented in \cite{jb1}. Roughly, a growth condition such as $(\log n)^{1+\delta}$ on $(\gamma_0\varepsilon_0/\alpha_0+\cdots+\gamma_{n-1}\varepsilon_{n-1}/\alpha_{n-1})$ is sufficient for transience, confirming the natural prevalence of transience results in the litterature on this model. 

\medskip
Extending \cite{jb1}, for the general model where $\alpha_n$ need not equal $\beta_n$, a full characterization of the recurrence regime was shown in \cite{jb2}. With some naturally generalized $\Phi(-m,n)$ (cf Definition \ref{maindef} below), the form of the criterion is the same, emphasizing the fact that the environment defines a new metrization of $\Z^{2}$. The model appears notably more general than the vertically flat one. Several examples were next given in \cite{jb2}. However, an empirical observation is that the methods employed for obtaining the structural results of \cite{jb1,jb2} are of very different nature than that used to treat examples. The analysis is in fact naturally divided in two parts, the second one never entering the mechanism of the random walk itself. The latter consists in studying fine properties of certain ergodic sums and is a source of interesting and difficult problems, for example closely related to temporal limit theorems and generalizations (cf Dolgopyat-Sarig \cite{dolgosarig}). In \cite{jb3}, for the general model, the particular case when the transition laws are independent was studied in detail, with a precise quantification of the non-surprising fact that the transience regime largely prevails in the set of parameters. 

\medskip
The purpose of the present article is to complete \cite{jb1,jb2,jb3}, by extending the applications of \cite{jb1,jb2}. We study for both the vertically flat and the general model the case when the transition laws are described by functions defined above an irrational rotation on the one-dimensional torus.

\section{Preliminaries}
Let $\T=\R\backslash\Z$ be the one-dimensional torus. Unless otherwise stated, functions are defined on $\T$, with arguments understood modulo one. We write $\Vert x\Vert$ for the distance of a real $x$ to $\Z$.

\medskip
We first recall classical facts about continued fractions. On this topic, one may consult Khinchin's book \cite{khinchin_book}. Any irrational $0<\theta<1$ admits an infinite continued fraction expansion:

$$\theta=\cfrac{1}{a_1+\cfrac{1}{a_2+\cdots}}=[0,a_1,a_2,\cdots],$$

\medskip
\noindent
where the partial quotients $(a_i)_{i\geq1}$ are integers $\geq1$, obtained by successive iterations of the Gauss map $x\longmapsto \{1/x\}$, starting from $\theta$. The convergents $(p_n/q_n)_{n\geq 1}$ of $\theta$ are the truncations $[0,a_1,a_2,\cdots,a_n]=p_n/q_n$, for $n\geq1$, of this continued fraction. The numerators $(p_n)$ and denominators $(q_n)$ check the same recurrence relation:

$$p_{n+1}=a_{n+1}p_n+p_{n-1},~~~q_{n+1}=a_{n+1}q_n+q_{n-1},~~n\geq0,$$  

\medskip
\noindent
with initial data $p_0=0,p_{-1}=1$ and $q_0=1,q_{-1}=0$. Classically (cf \cite{khinchin_book}, chap. 1):

\begin{equation}
\label{kh1}
\frac{1}{2q_{n+1}}\leq\frac{1}{q_n+q_{n+1}}\leq \Vert q_n\theta\Vert\leq\frac{1}{q_{n+1}}.
\end{equation}

\medskip
Fixing $\theta\not\in\Q$, we consider the rotation $Tx=x+\theta\mod1$ on $\T$ and write $T^nf$ for $f\circ T^n$, for any $f : \T\rightarrow\R$. We also use cocycle notations, for $x\in\T$:

$$f_n(x)=\left\{{\begin{array}{cc}
f(x)+\cdots+f(T^{n-1}x),&n\geq1,\\
0,&n=0,\\
-f(T^{n}x)-\cdots-f(T^{-1}x),&~n\leq-1.\end{array}
}\right.$$

\medskip
\noindent
An important property is that $f_{n+p}(x)=f_n(x)+T^nf_p(x)$, for any $x\in \T$ and $n,p\in\Z$.

\medskip
A function $f :\T\longmapsto\R$ with bounded variation will be said BV, with total variation written as $V(f)$. When $f$ is BV, with $\int_{\T}f(x)~dx=0$, the Denjoy-Koksma inequality says that:

\begin{equation}
\label{dk}
\vert f_{q_n}(x)\vert\leq V(f),~~~~~~~~~n\geq1,x\in\T.
\end{equation}

\medskip
Let us now recall known facts concerning Ostrowski's expansions (cf Beck \cite{beck}, p. 23). Every integer $q_m\leq n<q_{m+1}$ can be represented as:

\begin{equation}
\label{ostro}
n=\sum_{0\leq k\leq m}b_kq_k,
\end{equation}
 
\medskip
\noindent
with $0\leq b_0\leq a_1-1$, $0\leq b_j\leq a_{j+1}$, $1\leq j<m$, and $1\leq b_m\leq a_{m+1}$. Indeed, $n=b_mq_m+r$, for some $0\leq r<q_m$ and $1\leq b_m\leq a_m$. Iterating this process furnishes the decomposition \eqref{ostro}. Setting $A_{-1}=0$ and $A_k=\sum_{0\leq j\leq k}b_jq_j$, for $0\leq k<m$, by \eqref{ostro}, we have for any function $f$:

$$f_n(x)=\sum_{k=0}^mf_{b_kq_k}(x+A_{k-1}\theta).$$

\medskip
\noindent
When $f$ is BV and centered, the Denjoy-Koksma inequality \eqref{dk} furnishes the upper-bound:

\begin{equation}
\label{inegdk}
\vert f_{n}(x)\vert\leq \sum_{0\leq k\leq m}\Vert f_{q_k}\Vert_{\infty}b_k\leq V(f)\sum_{0\leq k\leq m}b_k,~~~~x\in\T.
\end{equation}

\medskip
Set $\N=\{0,1,\cdots\}$. For $g :\N\rightarrow\R_+$ increasing to $+\infty$ and $x\geq g(0)$, let $g^{-1}(x)$ be the unique integer $n\geq0$ such that $g(n)\leq x<g(n+1)$. By definition:

\begin{equation}
\label{sand}
g(g^{-1}(x))\leq x<g(g^{-1}(x)+1).
\end{equation}

\smallskip
\noindent
Also, $g^{-1}(g(n))=n$, for large $n\in\N$. Finally, for $f,g:\N\rightarrow\R_+$, we write $g\preceq f$ if there exists a constant $C>0$ so that $g(n)\leq Cf(n)$, for large $n\in\N$. We write $f\asymp g$ if $g\preceq f$ and $f\preceq g$.

\section{The quasi-periodic vertically flat model}
We first consider the vertically flat model, i.e. $\alpha_n=\beta_n=(1-\gamma_n)/2,~n\in\Z$. As a preliminary remark, we discuss the case when the sequence $(\varepsilon_n\gamma_n/(1-\gamma_n))_{n\in\Z}$ is periodic.

\begin{prop}

$ $

\noindent
For the vertically flat model, let $(\varepsilon_n\gamma_n/(1-\gamma_n))_{n\in\Z}$ be periodic with period $N\geq1$. Then the random walk is either recurrent or transient, according to whether:

$$\sum_{0\leq n<N}\varepsilon_n\gamma_n/(1-\gamma_n)=0\mbox{ or }\not=0.$$
\end{prop} 

\noindent
This follows from \cite{jb1}, respectively Prop 1.4 $i)$ and Corollary 1.3 $i)$. This extends the case of the Campanino-Petritis model in \cite{cp}, when $\mu_n=\delta_{x_n}$ with $x_n=\varepsilon_n=(-1)^n$, $\alpha_n=\beta_n=\gamma_n=1/3$. 

\medskip
Turning to quasi-periodic situations, we shall generalize \cite{jb1}, Prop. 1.5, giving in particular a better understanding of the Campanino-Petritis model in this quasi-periodic context.

\begin{thm}
\label{thm1} 

$ $

\noindent
Let $\theta=[0,a_1,a_2,\cdots]\not\in\Q$, with $\sum_{n\geq1}\log(1+a_n)/(a_1+\cdots+a_n)=+\infty$, and $Tx=x+\theta\mod1$ on $\T$. Let $f:\T\rightarrow\R$ be piecewise $K$-Lipschitz, with zero mean. Under Hypothesis \ref{hypo}, let $\alpha_n=\beta_n$ and $\varepsilon_n\gamma_n/(1-\gamma_n)=f(n\theta)$, $n\in\Z$. Then the random walk is recurrent.
\end{thm}

\noindent
\begin{remark}
This is shown in \cite{jb1}, Prop. 1.5, for $f=(1_{[0,1/2)}-1_{[1/2,1)})/2$, corresponding to the Campanino-Petritis model with $\mu_n=\delta_{x_n}$ and $x_n=1_{[0,1/2)}(n\theta)-1_{[1/2,1)}(n\theta)$. As a consequence of the theorem, the random walk is recurrent when $x_n=1_{[0,1/2)}(x+n\theta)-1_{[1/2,1)}(x+n\theta)$, $n\in\Z$, for any $x\in\T$, taking $f(.+x)$. As noticed in \cite{jb2}, Prop. 7.1, the condition on $\theta$ is generic in measure, since $\sum_{n\geq1}1/(a_1+\cdots+a_n)=+\infty$, for a.-e. $\theta$, cf Khinchin \cite{khinchin}. 
\end{remark}

\medskip
The other direction is in general more delicate, since requiring lower bounds on the ergodic sums. We just give an example.

\begin{prop}
\label{prop1}

$ $

\noindent
Let $\theta=[0,a_1,a_2,\cdots]\not\in\Q$, with $a_1$ odd and $a_n$ even for $n\geq2$. We suppose that for some $\delta>1$, $a_{n+1}\geq (a_n)^{\delta}$, for large $n$. Let $f=1_{[0,1/2)}-1_{[1/2,1)}$. Under Hypothesis \ref{hypo}, let $\alpha_n=\beta_n$ and $\varepsilon_n\gamma_n/(1-\gamma_n)=f(x+n\theta)$, $n\in\Z$, for some $x\in\T$. Then, for Lebesgue almost-every $x\in\T$, the random walk is transient.
\end{prop}

\noindent
\begin{remark}
In the last proposition, one may take for example the angle $\theta\in(0,1)$ defined by the partial quotients $a_n=2^{2^{n-1}-1}$, $n\geq1$. 
\end{remark}

\subsection{Proof of Theorem \ref{thm1}}

Fix $\theta=[0,a_1,a_2,\cdots,\cdots]\not\in\Q$, $Tx=x+\theta\mod1$ on $\T$ and $f$, as in the statement of the theorem. Using cocycle notations $(f_n(x))_{n\in\Z}$, introduce for $n\geq1$ the following positive functions $\varphi(n)$ and $\varphi_+(n)$ such that:

$$\varphi^2(n)=n^2+\sum_{-n\leq k<\ell\leq n}(f_{\ell}(x)-f_k(x))^2\mbox{ and }\varphi_+^2(n)=n^2+\sum_{-n\leq k< \ell\leq n,k\ell\geq0}(f_{\ell}(x)-f_k(x))^2.$$

\smallskip
\noindent
The dependence on $x$ of $\varphi(n)$ and $\varphi_+(n)$ is implicit. Obviously, $n\leq \varphi_+(n)\leq\varphi(n)$. The next lemma gives some control in the other direction.

\begin{lemme} 
\label{majo}

$ $

\noindent
There exists a constant $C_0>0$, uniform on $x\in\T$, such that for all $n\geq1$ and $1\leq m\leq 4a_{n+1}$:

$$\varphi^2(mq_n)\leq 2\varphi^2_+(mq_n)+C_0m^4q_n^2.$$
\end{lemme}

\noindent
{\it Proof of the lemma:}

\noindent
{\it Step 1.}  In the sequel, we simplify $f_n(x)$ into $f_n$, ${n\in\Z}$. Setting $A=\sum_{-n\leq k\leq -1,1\leq \ell\leq n}(f_{\ell}-f_k)^2$, we have $\varphi^2(n)=\varphi^2_+(n)+A$. Then:

\begin{eqnarray}
\label{edun}A&=&n\sum_{1\leq \ell\leq n}f_{\ell}^2+n\sum_{-n\leq k\leq -1}f_k^2-2\sum_{-n\leq k\leq -1}f_k\sum_{1\leq \ell\leq n}f_{\ell}.
\end{eqnarray}

\noindent
We next have:

\begin{equation}
\label{edun2}
-2\sum_{-n\leq k\leq -1}f_k\sum_{1\leq \ell\leq n}f_{\ell}=\left(\sum_{1\leq \ell\leq n}(f_{\ell}-f_{-\ell})\right)^2-\left(\sum_{1\leq \ell\leq n}f_{\ell}\right)^2-\left(\sum_{1\leq \ell\leq n}f_{-\ell}\right)^2.
\end{equation}

\smallskip
\noindent
Now, classically:

\begin{eqnarray}
\label{edun3}
\sum_{1\leq k<\ell\leq n}(f_{\ell}-f_k)^2&=&\sum_{2\leq \ell\leq n}({\ell}-1)f_{\ell}^2+\sum_{1\leq \ell\leq n-1}(n-\ell)f_{\ell}^2-2\sum_{1\leq k<\ell\leq n}f_kf_{\ell}\nonumber\\
&=&\sum_{1\leq \ell\leq n}(\ell-1)f_{\ell}^2+\sum_{1\leq \ell\leq n}(n-\ell)f_{\ell}^2-2\sum_{1\leq k<\ell\leq n}f_kf_{\ell}\nonumber\\
&=&n\sum_{1\leq \ell\leq n}f_{\ell}^2-\left(\sum_{1\leq \ell\leq n}f_{\ell}\right)^2.\end{eqnarray}

\smallskip
\noindent
Proceeding symmetrically for the other part of $A$, we obtain from \eqref{edun}, \eqref{edun2} and \eqref{edun3}:

$$A=\sum_{1\leq k<\ell\leq n}(f_{\ell}-f_k)^2+\sum_{-n\leq k<\ell\leq -1}(f_{\ell}-f_k)^2+\left(\sum_{1\leq \ell\leq n}(f_\ell-f_{-\ell})\right)^2.$$

\medskip
\noindent
Consequently:

$$\varphi^2(n)\leq 2\varphi^2_+(n)+\left(\sum_{1\leq \ell\leq n}(f_{\ell}-f_{-\ell})\right)^2.$$

\medskip
\noindent
{\it Step 2.} Let $n\geq 1$ and $1\leq m\leq 4a_{n+1}$. Setting $B=\sum_{1\leq \ell\leq mq_n}(f_{\ell}-f_{-\ell})$, we have:

\begin{eqnarray}
B&=&\sum_{0\leq u<m}\sum_{1\leq \ell\leq q_n}(f_{uq_n+\ell}-f_{-uq_n-\ell})=\sum_{0\leq u<m}\sum_{1\leq \ell\leq q_n}(f_{uq_n}-f_{-uq_n}+T^{uq_n}f_{\ell}-T^{-uq_n}f_{-\ell}).\nonumber
\end{eqnarray}

\medskip
\noindent 
Using Denjoy-Koksma's inequality \eqref{dk}, for any integer $0\leq u<m$, we have $\vert f_{uq_n}(x)\vert\leq u V(f)$, idem for $f_{-uq_n}(x)$. As a result:

\begin{eqnarray}
B&=&O(m^2q_n)+\sum_{0\leq u<m}\sum_{1\leq \ell\leq q_n}(T^{uq_n}f_{\ell}-T^{-uq_n}f_{-\ell}).\nonumber
\end{eqnarray}

\medskip
\noindent 
Fixing any integer $0\leq u<m$, we have :

\begin{eqnarray}
\sum_{1\leq \ell\leq q_n}(T^{uq_n}f_{\ell}-T^{-uq_n}f_{-\ell})&=&\sum_{k=0}^{q_n-1}(q_n-k)f(x+uq_n\theta+k\theta)+\sum_{k=1}^{q_n}(q_n+1-k)f(x-uq_n\theta-k\theta).\nonumber
\end{eqnarray}

\medskip
\noindent
Using the Denjoy-Koksma inequality \eqref{dk} for the $(q_n+1)$-term in the second sum on the right hand side and making a change of variable in the first one, we get:

\begin{eqnarray}\sum_{1\leq \ell\leq q_n}(T^{uq_n}f_{\ell}-T^{-uq_n}f_{-\ell})&=&\sum_{k=1}^{q_n}k(f(x+uq_n\theta+(q_n-k)\theta)-f(x-uq_n\theta-k\theta))+O(q_n)\nonumber\\
&=&\sum_{k=1}^{q_n}k(f(x_+^u-k\theta)-f(x_-^u-k\theta))+O(q_n),\nonumber\end{eqnarray}

\medskip
\noindent
when setting $x^u_{+}=x+(u+1)q_n\theta$ and $x^u_{-}=x-uq_n\theta$. By \eqref{kh1} and the hypothesis $m\leq 4a_{n+1}$:

$$\|x^u_+-x^u_-\|\leq (2u+1)\Vert q_n\theta\Vert\leq (8a_{n+1}+1)/q_{n+1}\leq 9/q_n.$$

\smallskip
\noindent
Denote by $[x^u_-,x_+^u]$ the short interval on $\T$ between $x_-^u$ and $x_+^u$ and by $D$ the number of discontinuities of $f$. Recall also that there is exactly one $k\theta\mod1$, $1\leq k\leq q_n$, in each interval $[\ell/q_n,(\ell+1)/q_n)$, $0\leq\ell<q_n$, on $\T$. As a result, for a given discontinuity of $f$, there at most $10$ values of $1\leq k\leq q_n$ such that $[x^u_-,x_+^u]-k\theta$ contains this discontinuity. Hence, using that $f$ is $K$-Lipschitz on intervals containing no discontinuity, an upper-bound for the last sum is:

$$D\times 10q_n\times 2\Vert f\Vert_{\infty}+q_n^2\times K\times\frac{9}{q_n}=O(q_n),$$

\medskip
\noindent
As a result, $B=O(m^2q_n)+O(mq_n)=O(m^2q_n)$, which ends the proof of the lemma.

\fin

\bigskip
We turn to the proof of Theorem \ref{thm1}. We assume that $\alpha_n=\beta_n$ and that $\mu_n$ and thus $\varepsilon_n$ are such that $\varepsilon_n\gamma_n/(1-\gamma_n)=f(x+n\theta)$, $n\in\Z$, for some $x\in\T$. The statement of the theorem corresponds to $x=0$. Introduce the following definition, due to Feller (1969):

\begin{defi}
\label{domiv}

$ $

\noindent
A non-decreasing function $g:\R_+\rightarrow\R_+$ satisfies dominated variation, if there exists a constant $C>0$ so that for large $x>0$,~$g(2x)\leq Cg(x)$. Hence, iterating, for all $K>0$, there exists $C_K$ so that for large $x>0$, $g(Kx)\leq C_Kg(x)$.
 
\end{defi}

\noindent
In \cite{jb1}, setting $R_k^{\ell}=\sum_{k\leq i\leq \ell}\varepsilon_i\gamma_i/(1-\gamma_i)$, for integers $k\leq \ell$, we considered the two functions:

$$\Phi^2(n)=n^2+\sum_{-n\leq k\leq \ell\leq n}(R_k^{\ell})^2,~\Phi_+^2(n)=n^2+\sum_{-n\leq k\leq \ell\leq n,k\ell>0}(R_k^{\ell})^2,$$ 

\noindent
The following results were then established:

\begin{thm}(\cite{jb1}, Lemma 6.1, Theorem 1.2. and Corollary 1.3 $i)$)

\noindent
1) The functions $\Phi$ and $\Phi_+$ satisfy dominated variation with a constant $C=C(\eta)$ depending only on $\eta$, where $\eta$ is introduced in Hypothesis \ref{hypo}.

\medskip
\noindent
2) The random walk is recurrent if and only if $\sum_{n\geq1}n^{-2}(\Phi^{-1}(n))^2/\Phi_+^{-1}(n)=+\infty$.

\medskip
\noindent
3) The condition $\sum_{n\geq1}1/\Phi_+(n)<+\infty$ is sufficient for transience.

\end{thm}

\noindent
Using that $f$ is bounded and that $\left|{|R_k^{\ell}|-|f_{\ell}(x)-f_k(x)|}\right|\leq \|f\|_{\infty}$, $k\leq \ell$, it is immediate that for some constant $C>0$ depending only on $\eta$ (hence uniform on $x\in\T$), for all $n\geq1$:

$$\Phi(n)/C\leq\varphi(n)\leq C\Phi(n)\mbox{ and }\Phi_+(n)/C\leq\varphi_+(n)\leq C\Phi_+(n).$$

\begin{cor} 
\label{corro}
$ $

\noindent
1) The functions $\varphi^{-1}$ and $\varphi^{-1}_+$ satisfy dominated variation, i.e. for any $K>0$, there exists a constant $C_K>0$, independent of $x\in\T$, so that for large $y>0$: 

\begin{equation}
\label{domvar}
\varphi^{-1}(Ky)\leq C_K\varphi^{-1}(y)\mbox{ and }\varphi_+^{-1}(Ky)\leq C_K\varphi_+^{-1}(y).
\end{equation}

\smallskip
\noindent
2) The random walk is recurrent if and only if $\sum_{n\geq1}n^{-2}(\varphi^{-1}(n))^2/(\varphi_+^{-1}(n))=+\infty$.

\medskip
\noindent
3) The condition $\sum_{n\geq1}1/\varphi_+(n)<+\infty$ is sufficient for transience.
\end{cor}

Let us reprove a concrete version of dominated variation of $\varphi^{-1}$ and $\varphi_+^{-1}$ in the following lemma. Concerning for example $\varphi_+$, we essentially show that $n \longmapsto \varphi^2_+(n)/n$, $n>0$, is non-decreasing. For $a<b$ in $\Z$ and $x\in\T$, let:

\begin{equation}
\label{psi}
\psi(a,b)=\sum_{a\leq k<\ell\leq b}(f_{\ell}(x)-f_k(x))^2,
\end{equation}

\smallskip
\noindent
where the dependence on $x$ is implicit on the left hand side.

\begin{lemme}
\label{triinv}

$ $

\noindent
Let integers $a<b<c$ and $x\in\T$. Then:

\begin{eqnarray}
\frac{\psi(a,c)}{c-a}\geq \frac{\psi(a,b-1)}{b-a}+\frac{\psi(b+1,c)}{c-b}.\nonumber
\end{eqnarray}

\noindent
Also, for large $n$, uniformly in $x\in\T$:

\begin{equation}
\label{limin}
\varphi_+(2n)\geq 2^{1/4}\varphi_+(n).
\end{equation}
\end{lemme}

\noindent
{\it Proof of the lemma:}

\noindent
We write $f_k$ in place of $f_k(x)$. Decompose $\psi(a,c)=\psi(a,b)+\psi(b,c)+\sum_{a\leq k< b< \ell\leq c}(f_{\ell}-f_k)^2$ and then expend:

\begin{eqnarray}
\sum_{a\leq k< b< \ell\leq c}(f_{\ell}-f_k)^2&=&(b-a)\sum_{b<\ell\leq c}f_{\ell}^2+(c-b)\sum_{a\leq k< b}f_k^2-2\sum_{b<\ell\leq c}f_{\ell}\sum_{a\leq k< b}f_k.\nonumber
\end{eqnarray}

\medskip
\noindent
As before, $\psi(b+1,c)=(c-b)\sum_{b< \ell\leq c}f_{\ell}^2-(\sum_{b< \ell\leq c}f_{\ell})^2$ and $\psi(a,b-1)=(b-a)\sum_{a\leq k< b}f_k^2-(\sum_{a\leq k< b}f_k)^2$. When substituting:

\begin{eqnarray}
\sum_{a\leq k< b< \ell\leq c}(f_{\ell}-f_k)^2&=&\frac{b-a}{c-b}\left({\psi(b+1,c)+(\sum_{b< \ell\leq c}f_{\ell})^2}\right)\nonumber\\
&+&\frac{c-b}{b-a}\({\psi(a,b-1)+(\sum_{a\leq k< b}f_k)^2}\)-2\sum_{b<\ell\leq c}f_{\ell}\sum_{a\leq k< b}f_k.\nonumber
\end{eqnarray}

\medskip
\noindent
As a consequence:

\begin{eqnarray}
\label{ineg2}
\psi(a,c)&\geq&\frac{c-a}{c-b}\psi(b+1,c)+\frac{c-a}{b-a}\psi(a,b-1)\nonumber\\
&+&\left({\sqrt{\frac{b-a}{c-b}}\sum_{b< \ell\leq c}f_{\ell}- \sqrt{\frac{c-b}{b-a}}\sum_{a\leq k< b}f_k}\right)^2,
\end{eqnarray}

\medskip
\noindent
giving the result. Concerning the last part of the lemma, for any integer $n\geq1$, we have:

$$\psi(0,2n)\geq \frac{2n}{n+1}\psi(0,n),~\psi(-2n,0)\geq\frac{2n}{n+1}\psi(-n,0).$$

\smallskip
\noindent
Since $\varphi_+^2(n)=n^2+\psi(-n,0)+\psi(0,n)$, we get \eqref{limin}. 

\fin

\bigskip
\noindent
\begin{remark}
As a complement, let us observe that:

\begin{equation}
\label{ineg1}
\frac{c-a}{c-b}\psi(b+1,c)+\frac{c-a}{b-a}\psi(a,b-1)\geq (\sqrt{\psi(a,b-1)}+\sqrt{\psi(b+1,c)})^2,
\end{equation}

\medskip
\noindent
as this is equivalent to the true relation:

$$\frac{c-b}{b-a}\psi(a,b-1)+\frac{b-a}{c-b}\psi(b+1,c)\geq2\sqrt{\psi(a,b-1)}\sqrt{\psi(b+1,c)}.$$
\end{remark}

\medskip
\noindent
Using \eqref{ineg1} in \eqref{ineg2}, this thus implies some reverse triangular inequality:

$$\sqrt{\psi(a,c)}\geq\sqrt{\psi(a,b-1)}+\sqrt{\psi(b+1,c)}.$$

\medskip
We start the proof of Theorem \ref{thm1}. A corollary of Lemma \ref{majo} is that there exists a constant $C_0>0$, independent of $x\in\T$, such that for all $n\geq1$ and $1\leq m\leq 4a_{n+1}$:

\begin{equation}
\label{corot}
\varphi(mq_n)\leq C_0(\varphi_+(mq_n)+m^2q_n).
\end{equation}

\bigskip
Let now $n\geq 1$ and $\ell\geq0$ be such that $2^{\ell}\leq 4a_{n+1}$. We make the following discussion: 

\medskip
\noindent
- {\it Case 1.} If $\varphi_+(2^{\ell}q_n)\geq 2^{2\ell}q_n$, then, using \eqref{corot} and next \eqref{domvar} at the end:

\begin{eqnarray}
\varphi^{-1}(\varphi_+(2^{\ell}q_n))\geq \varphi^{-1}((\varphi_+(2^{\ell}q_n)+2^{2\ell}q_n)/2)&\geq& \varphi^{-1}(\varphi(2^{\ell}q_n)/(2C_0))\nonumber\\
&\geq&\frac{\varphi^{-1}(\varphi(2^{\ell}q_n))}{C_{2C_0}}=\frac{2^{\ell}q_n}{C_{2C_0}}.\nonumber
\end{eqnarray}

\medskip
\noindent
By Ostrowski's expansion \eqref{ostro}, for $|k|\leq 2^{\ell} q_n\leq 4q_{n+1}$, $|f_k|\leq 4\times V(f)(a_1+\cdots+a_{n+1})$. Hence, for some $C>0$, $\varphi_+(2^{\ell}q_n)\leq C2^{\ell}q_n(a_1+\cdots+a_{n+1})$. Thus, with $C_1=1/(CC_{2C_0})$, for large $n$:

\begin{equation}
\label{c1}
\frac{(\varphi^{-1}(\varphi_+(2^{\ell}q_n)))^2}{2^{\ell}q_n\varphi_+(2^{\ell}q_n)}\geq\frac{1}{C_{2C_0}} \frac{2^{\ell}q_n}{\varphi_+(2^{\ell}q_n)}\geq \frac{C_1}{a_1+\cdots+a_{n+1}}.
\end{equation}

\medskip
\noindent
- {\it Case 2.} Suppose that $\varphi_+(2^{\ell}q_n)< 2^{2\ell}q_n$. As $\varphi_+(2^{\ell}q_n)\geq 2^{\ell}q_n$, there exists $0\leq \ell'\leq \ell$ so that $2^{2\ell'}q_n\leq \varphi_+(2^{\ell}q_n)<2^{2(\ell'+1)}q_n$. We get in this case, using \eqref{corot} and \eqref{domvar}:

\begin{eqnarray}
\varphi^{-1}(\varphi_+(2^{\ell}q_n))&\geq& \varphi^{-1}(2^{2\ell'}q_n)\geq\varphi^{-1}((\varphi_+(2^{\ell'}q_n)+2^{2\ell'}q_n)/5)
\geq \varphi^{-1}(\varphi(2^{\ell'}q_n)/(5C_0))\nonumber\\
&\geq&\frac{\varphi^{-1}(\varphi(2^{\ell'}q_n))}{C_{5C_0}}=\frac{2^{\ell'}q_n}{C_{5C_0}}.\nonumber
\end{eqnarray}

\medskip
\noindent
As a result, we can write in this case, redefining $C_1:=\min\{C_1,1/(16C_{5C_0})\}>0$:

$$\frac{(\varphi^{-1}(\varphi_+(2^{\ell}q_n)))^2}{2^{\ell}q_n\varphi_+(2^{\ell}q_n)}\geq \frac{1}{C_{5C_0}}\frac{2^{2\ell'}q^2_n}{2^{\ell}q_n\varphi_+(2^{\ell}q_n)}\geq \frac{1}{4C_{5C_0}2^{\ell}}\geq \frac{C_1}{a_1+\cdots+a_{n+1}}.$$

\medskip
Hence, \eqref{c1} is true for any $\ell\geq 0$ such that $2^{\ell}\leq 4a_{n+1}$. This now gives, for large $n>0$:

\begin{eqnarray}
\sum_{\varphi_+(q_n)\leq k<\varphi_+(4a_{n+1}q_n)}\frac{1}{k^2}\frac{(\varphi^{-1}(k))^2}{\varphi_+^{-1}(k)}&\geq&\sum_{0\leq \ell\leq 1+\log_2a_{n+1}}\sum_{\varphi_+(2^{\ell}q_n)\leq k<\varphi_+(2^{\ell+1}q_{n})}\frac{1}{k^2}\frac{(\varphi^{-1}(k))^2}{\varphi_+^{-1}(k)}\nonumber\\
&\geq&\sum_{0\leq \ell\leq1+\log_2a_{n+1}}\frac{(\varphi^{-1}(\varphi_+(2^{\ell}q_n)))^2}{\varphi_+^{-1}(\varphi_+(2^{\ell+1}q_n))}\sum_{\varphi_+(2^{\ell}q_n)\leq k<\varphi_+(2^{\ell+1}q_{n})}\frac{1}{k^2}\nonumber\\
&\geq&\sum_{0\leq \ell\leq 1+\log_2a_{n+1}}\frac{(\varphi^{-1}(\varphi_+(2^{\ell}q_n)))^2}{2^{\ell+1}q_n}\sum_{\varphi_+(2^{\ell}q_n)\leq k<\varphi_+(2^{\ell+1}q_{n})}\frac{1}{k^2}.\nonumber
\end{eqnarray}

\medskip
\noindent
Using now relation \eqref{c1}, we arrive at:

\begin{eqnarray}
\sum_{\varphi_+(q_n)\leq k<\varphi_+(4a_{n+1}q_n)}\frac{1}{k^2}\frac{(\varphi^{-1}(k))^2}{\varphi_+^{-1}(k)}
&\geq&\frac{C_1}{2(a_1+\cdots+a_{n+1})}\sum_{0\leq \ell\leq 1+\log_2a_{n+1}}\varphi_+(2^{\ell}q_n)\sum_{\varphi_+(2^{\ell}q_n)\leq k<\varphi_+(2^{\ell+1}q_{n})}\frac{1}{k^2}.\nonumber\end{eqnarray}

\medskip
\noindent
Using that $1/k^2\geq1/k-1/(k+1)$, we obtain:

$$\sum_{\varphi_+(2^{\ell}q_n)\leq k<\varphi_+(2^{\ell+1}q_{n})}\frac{1}{k^2}\geq\frac{1}{\varphi_+(2^{\ell}q_n)+1}-\frac{1}{\varphi_+(2^{\ell+1}q_n)-1}.$$

\medskip
\noindent
Now, for large $n$ (uniformly in $\ell\geq 0$), applying \eqref{limin}, we have $\varphi_+(2^{\ell+1}q_n)\geq 2^{1/4}\varphi_+(2^{\ell}q_n)$. As a result, for large $n$, we obtain:

$$\sum_{\varphi_+(2^{\ell}q_n)\leq k<\varphi_+(2^{\ell+1}q_{n})}\frac{1}{k^2}\geq\frac{1-2^{-1/4}}{2}\frac{1}{\varphi_+(2^{\ell}q_n)}.$$

\smallskip
\noindent
This thus furnishes, for large $n>0$:

\begin{eqnarray}
\sum_{\varphi_+(q_n)\leq k<\varphi_+(4a_{n+1}q_n)}\frac{1}{k^2}\frac{(\varphi^{-1}(k))^2}{\varphi_+^{-1}(k)}
&\geq&\frac{C_1(1-2^{-1/4})}{4(a_1+\cdots+a_{n+1})}(1+\log_2 a_{n+1}).\nonumber\end{eqnarray}

\medskip
\noindent
Finally, notice that $4a_{n+1}q_n\leq q_{n+5}$, for $n\geq0$. Hence, for large $N>0$:

$$\sum_{k\geq \varphi_+(q_N)}\frac{(\varphi^{-1}(k))^2}{\varphi_+^{-1}(k)}\geq\frac{1}{5}\sum_{n\geq N}\sum_{\varphi_+(q_n)\leq k<\varphi_+(q_{n+5})}\frac{(\varphi^{-1}(k))^2}{\varphi_+^{-1}(k)}\geq\frac{C_1(1-2^{-1/4})}{5}\sum_{n\geq N}\frac{(1+\log_2 a_{n+1})}{4(a_1+\cdots+a_{n+1})}.$$

\smallskip
\noindent
By hypothesis, the series on the right hand side diverges. From Corollary \ref{corro}, we conclude that the random walk is recurrent.

\fin

\subsection{Proof of Proposition \ref{prop1}} 

Let us place in the context of this proposition, namely $\alpha_n=\beta_n$ and $\varepsilon_n\gamma_n/(1-\gamma_n)=f(x+n\theta)$, $n\in\Z$, for some $x\in\T$. Here $f=1_{[0,1/2)}-1_{[1/2,1)}$ and the angle $\theta=[0,a_1,a_2,\cdots]$ is such that $a_1$ is odd and $a_n$ is even for $n\geq2$, together with $a_{n+1}\geq (a_n)^{\delta}$, for large $n$, for some fixed $\delta>1$.

\medskip
From the relation $q_{n+1}=a_{n+1}q_n+q_{n-1}$, $n\geq0$, and $q_0=1$, $q_{-1}=0$, we recursively obtain that $q_n$ is odd, $n\geq1$. This implies that for all $n\geq1$ and $x\in\T$, $\vert f_{q_n}(x)\vert\geq1$, since the number of $(T^kx)_{0\leq k<q_n}$ that fall in the intervals $[0,1/2)$ and $[1/2,1)$ are different.

\medskip
Let $2/3<\beta<1$ and define $m_k=(a_{k+1})^{\beta}$, $k\geq1$. Introduce :

$$A_k=\{x\in\T,f_{mq_k}(x)=mf_{q_k}(x),~ 0\leq m\leq m_k\}.$$

\medskip
\noindent
Write ${\cal L}_{\T}$ for Lebesgue measure on $\T$ and denote by $[x,x+q_k\theta]$ the short interval on $\T$ determined by $x$ and $x+q_k\theta$ on $\T$. If $T^r[x,x+q_k\theta]$ does not contain either $0$ or $1/2$, for $0\leq r<m_kq_k$, then $x\in A_k$, because $f(T^rx)=f(T^{r+q_k}x)$, for $0\leq r<m_kq_k$. Hence, if $x\in\T\backslash A_k$, there exist $0\leq r<m_kq_k$ such that either $0$ or $1/2$ belongs to $T^r[x,x+q_k\theta]$. As a result, $\T\backslash  A_k\subset\cup_{0\leq r<m_kq_k,y\in\{0,1/2\}}([y-q_k\theta,y]-r\theta)$. Using \eqref{kh1}:

$${\cal L}_{\T}(\T\backslash A_k)\leq \frac{2m_kq_k}{q_{k+1}}\leq \frac{2m_k}{a_{k+1}}\leq 2(a_{k+1})^{-(1-\beta)}.$$

\smallskip
\noindent
Since $(a_k)$ grows at least geometrically, $\sum_ka_k^{-(1-\beta)}<+\infty$. By the first lemma of Borel-Cantelli, we deduce that for Lebesgue a.-e. $x$, $x\in A_k$ for large $k$.

\medskip
Let $N_k=a_1+\cdots+a_k$ and observe that $N_k\sim_{k\rightarrow+\infty}a_k$, since for large $k$, $a_k\geq (a_{k-1})^{\delta}$, with $\delta>1$. As $m_k=(a_{k+1})^{\beta}$, we now choose $\beta<1$ close enough to 1, so that $m_k\geq1000N_k$, for large $k$. Let now $m_k\geq m\geq 100N_k$. For $0\leq \ell,\ell'\leq mq_k$, we make the Euclidean divisions of $\ell,\ell'$ by $q_k$ : $\ell=aq_k+b$ and $\ell'=a'q_k+b'$, $0\leq b,b'<q_k$ and where $0\leq a,a'\leq m_k$. Then, almost-surely for large $k$, since $x\in A_k$, we have $f_{aq_k}(x)=af_{q_k}(x)$, $f_{a'q_k}(x)=a'f_{q_k}(x)$. Thus:

$$f_{\ell}(x)-f_{\ell'}(x)=(a-a')f_{q_k}(x)+T^{aq_k}f_b(x)-T^{a'q_k}f_{b'}(x).$$

\medskip
\noindent
Using the upper-bound \eqref{inegdk}, coming from Ostrowski's expansion \eqref{ostro}, for $\|f_b\|_{\infty}$ and $\|f_{b'}\|_{\infty}$, we have (since $V(f)=2$), $|T^{aq_k}f_b(x)|\leq 2 N_k$ and $|T^{a'q_k}f_{b'}(x)|\leq 2N_k$. Hence, a.-e., for large $k$, using the fact that $|f_{q_k}(x)|\geq1$:

\begin{equation}
\label{eucl}\vert f_{\ell}(x)-f_{\ell'}(x)\vert \geq \vert a-a'\vert \vert f_{q_k}(x)\vert -\vert T^{aq_k}f_b(x)-T^{a'q_k}f_{b'}(x)\vert \geq \vert a-a'\vert - 4N_k.
\end{equation}

\medskip
\noindent
Consequently, for $m_k\geq m\geq 100N_k$, by \eqref{eucl}:

\begin{eqnarray}
\label{tutur}
\varphi^2_+(mq_k)\geq\sum_{0\leq \ell'\leq \ell\leq mq_k}(f_{\ell}-f_{\ell'})^2&\geq&\sum_{0\leq a'<m/4; m/2<a<m;0\leq b,b'<q_k}(f_{aq_k+b}(x)-f_{a'q_k+b'}(x))^2\nonumber\\
&\geq&\sum_{0\leq a'<m/4; m/2<a<m;0\leq b,b'<q_k}(m/4-4N_k)^2\nonumber\\&\geq&(m/4)^2q_k^2(m/5)^2\geq m^4q_k^2/400.\end{eqnarray}

\medskip
In order to conclude the argument we shall apply Corollary \ref{corro}, 3), and show the convergence of $\sum_{n\geq1}1/\varphi_+(n)$. Let us write for a.-e. $x$ and large $K>0$: 

\begin{eqnarray}
\sum_{n\geq_K}\frac{1}{\varphi_+(n)}&\leq&\sum_{k\geq K}\sum_{1\leq m\leq a_{k+1} }\sum_{mq_k\leq n<(m+1)q_{k}}\frac{1}{\varphi_+(n)}\leq\sum_{k\geq K}\sum_{1\leq m\leq a_{k+1}}\frac{q_k}{\varphi_+(mq_k)}\nonumber\\
&\leq&\sum_{k\geq K}\[{\sum_{1\leq m\leq 100N_k}+\sum_{100N_k< m\leq m_{k}}+\sum_{m_k<m\leq a_{k+1}}}\]\({\frac{q_k}{\varphi_+(mq_k)}}\)\nonumber\\
&=&\sum_{k\geq K}~\[{~~~~~~U_k~~~~+~~~~~~~V_k~~~~~~+~~~~~W_k~~~~~}\].\nonumber\end{eqnarray}

\medskip
\noindent
1) Considering $\sum_{k\geq K}U_k$, we fix $k\geq K$ and $1\leq m\leq 100N_k$. Using first the function $\psi$ introduced in \eqref{psi}, before Lemma \ref{triinv}, we have:

\begin{equation}
\label{deco}
\varphi_+^2(mq_k)=(mq_k)^2+\psi(-mq_k,0)+\psi(0,mq_k).
\end{equation}

\smallskip
\noindent
By Lemma \ref{triinv} and next \eqref{tutur}, for some (next generic) constant $c>0$, a.-e. for large $k>0$:

$$\varphi_+(mq_k)\geq c\sqrt{mq_k/(m_{k-1}q_{k-1})}\varphi_+(m_{k-1}q_{k-1})\geq c\sqrt{mq_k/(m_{k-1}q_{k-1})}m^2_{k-1}q_{k-1}.$$

\medskip
\noindent
As a result:

\begin{equation}
\label{uk}
\varphi_+(mq_k)\geq c\sqrt{ma_k}m_{k-1}^{3/2}q_{k-1}\geq c\sqrt{m}a_k^{1/2+3\beta/2}q_{k-1}\geq\({ c\sqrt{m}a_k^{3\beta/2-1/2}}\)q_{k}.
\end{equation}

\medskip
\noindent
We obtain, a.-e., for large $K>0$, via \eqref{uk}, still for a generic $c>0$, using at the end that $N_k\sim_{k\rightarrow+\infty} a_k$, $\beta>2/3$ and that $(a_k)$ grows at least geometrically:

\begin{eqnarray}
\sum_{k\geq K}U_k&\leq& c\sum_{k\geq K}\sum_{1\leq m\leq 100N_k}\frac{1}{\sqrt{m}a_k^{3\beta/2-1/2}}\nonumber\\
&\leq& c\sum_{k\geq K}\frac{\sqrt{N_k}}{a_k^{3\beta/2-1/2}}
\leq c\sum_{k\geq K}\frac{\sqrt{a_k}}{a_k^{3\beta/2-1/2}}\leq c\sum_{k\geq1}a_k^{1-3\beta/2}<+\infty.\nonumber
\end{eqnarray}

\medskip
\noindent

\medskip
\noindent
2) Considering now $\sum_{k\geq K}V_k$, let $k\geq K$ and $100N_k<m\leq m_k$. Using \eqref{tutur}, we have a.-e. for large $K>0$ that $\varphi_+(mq_k)\geq m^2q_k/20$. So, for a generic constant $c>0$:

$$\sum_{k\geq K}V_k\leq c\sum_{k\geq K}\sum_{100N_k< m\leq m_{k}}\frac{1}{m^2}\leq c\sum_{k\geq K}\frac{1}{N_k}\leq c\sum_{k\geq K}a_k^{-1}<+\infty.$$

\medskip
\noindent
3) For $\sum_{k\geq K}W_k$, let $k\geq K$ and $m_k<m\leq a_{k+1}$. Again, using \eqref{deco}, Lemma \ref{triinv} and finally \eqref{tutur}, for $m_k$, we obtain, a.-e., for large $K>0$, for some generic constant $c>0$:

$$\varphi_+(mq_k)\geq  c\sqrt{m/m_k}\varphi_+(m_kq_k)\geq  c\sqrt{m/m_k}m^2_kq_k\geq c\sqrt{m}m_k^{3/2}q_k.$$ 

\smallskip
\noindent
As a consequence, we can write, for some generic $c>0$, using at the end that $\beta>2/3>1/3$:

$$\sum_{k\geq K}W_k\leq c\sum_{k\geq K}\sum_{m_k< m\leq a_{k+1}}\frac{1}{\sqrt{m}m_k^{3/2}}\leq c\sum_{k\geq K}\frac{\sqrt{a_{k+1}}}{m_k^{3/2}}\leq c\sum_{k\geq K}a_{k+1}^{1/2-3\beta/2}<+\infty.$$

\medskip
\noindent
This completes the proof of the proposition.

\fin

\section{The quasi-periodic general case}

Let again $Tx=x+\theta\mod 1$ be an irrational rotation on $\T$. The basic assumption in this section will be that for some BV functions $f,g :\T\rightarrow\R$, with $f$ centered, for some $x\in\T$:

$$\beta_n/\alpha_n=e^{f(T^{n-1}x)}, \gamma_n\varepsilon_n/\alpha_n=g(T^nx),~n\in\Z.$$

\medskip
\noindent
Implicitly, Hypothesis \ref{hypo} will always be realized, uniformly in $x\in\T$. Introduce some definitions.  

\begin{defi} 
\label{4.1}

$ $

\noindent
Fixing $x\in\T$, let $\rho_n=\rho_n(x)=e^{f_n(x)}$, $n\in\Z$. In other words:

$$\rho_n=\left\{{\begin{array}{cc}\frac{\beta_1}{\alpha_1}\cdots \frac{\beta_n}{\alpha_n},&~n\geq1,\\
1,&~n=0,\\
\frac{\alpha_{n+1}}{\beta_{n+1}}\cdots\frac{\alpha_0}{\beta_0},&~n\leq -1.\\
\end{array}}\right.$$

\smallskip
\noindent
For $n\geq0$, let:

$$v_+(n)=\sum_{0\leq k\leq n}\rho_k\mbox{ and }v_-(n)=(q_0/p_0)\sum_{-n-1\leq k\leq -1}\rho_k.$$

\smallskip
\noindent
In the same way, introduce for $n\geq0$:

$$w_+(n)=\sum_{0\leq k\leq n}1/\rho_k\mbox{ and }w_-(n)=(p_0/q_0)\sum_{-n-1\leq k\leq -1}1/\rho_k.$$

\end{defi}

\medskip
As already mentioned in the Introduction, for the random walk to be recurrent, the vertical random walk has first to be. The necessary and sufficient condition for this (cf \cite{KMG}) is:

$$\lim_{n\rightarrow+\infty} v_+(n)=+\infty~~~\mbox{ and    }~~\lim_{n\rightarrow+\infty} v_-(n)=+\infty.$$

\medskip
\noindent
Since $f$ is BV and centered, the Denjoy-Koksma inequality \eqref{dk} says that $\vert f_{\pm q_n}(x)\vert\leq V(f)$, for any $x\in\T$. As a result, $\rho_n=e^{f_n(x)}$ does not go to zero, neither as $n\rightarrow+\infty$, nor as $n\rightarrow-\infty$, implying that the two previous conditions hold. 

\bigskip
Some quasi-invariant measures on $\T$ will play a role. In the sequel, we consider the space of Borel probability measures on $\T$, equipped with its usual (metrizable) weak-$*$ topology (using continuous $w:\T\rightarrow\R$, as test functions), for which this space is compact. We denote by $T\nu$ the image by $T$ of a Borel probability measure $\nu$ on $\T$. By definition, $\int w~dT\nu=\int Tw~d\nu$, for any bounded measurable $w: \T\rightarrow\R$. Let us recall the following folklore result. 

\begin{thm} (\cite{cg}, Prop. 5.8, or \cite{anss}, Prop. 1.1.) 
\label{cgg}

\noindent
Let $h : \T\rightarrow\R$ be BV and centered. There exists a unique Borel probability measure $\nu_h$ on $\T$ such that $dT\nu_h=e^{T^{-1}h}d\nu_h$. This measure has no atom.
\end{thm}

\noindent
The proof of existence, relying on non-atomicity, is incomplete in \cite{cg} and too abstract in \cite{anss}. We choose to reprove existence and atomicity in an elementary way below. 

\medskip
\noindent
Also, it is well-known (cf \cite{cg}) that $\nu_h$ is absolutely continuous with respect to Lebesgue measure ${\cal L}_{\T}$ if and only if $h=\log u-\log Tu$, for some ${\cal L}_{\T}$-integrable $u>0$, otherwise it is singular. When $\nu_h$ is as in Theorem \ref{cgg}, notice the following relation, for any bounded measurable $w:\T\rightarrow\R$:

\begin{equation} 
\label{inva}
\int_{\T} w~d\nu_h=\int_{\T} e^{-h}Tw~d\nu_h.
\end{equation}

\medskip
We shall show in this section:

\begin{thm}
\label{thm2} 

$ $

\noindent
Let $\theta\not\in\Q$  and $Tx=x+\theta\mod1$ on $\T$. Let BV functions $f,g :\T\rightarrow\R$, with $f$ centered. Suppose that $\beta_n/\alpha_n=e^{f(T^{n-1}x)}$ and $\gamma_n\varepsilon_n/\alpha_n=g(T^nx), n\in\Z$, for some $x\in\T$.

\medskip
\noindent
i) Suppose that $\int_{\T}gd\nu_{f}\not=0$. Then for all $x\in\T$, the random walk is transient.

\medskip
\noindent
ii) Assume that $g=h-e^{-f}Th$, for some bounded $h$. Introduce:

\smallskip
1) $f=u-Tu$, with $e^{u}\in L^1({\cal L}_{\T})$. 

2) $f(x_0+x)=f(x_0-x)$, for some $x_0\in\T$ and ${\cal L}_{\T}$-a.-e. $x\in\T$. 

\smallskip
\noindent
Then, under either condition 1) or 2), for ${\cal L}_{\T}$-a.-e. $x$, the random walk is recurrent.

\end{thm} 

\noindent
\begin{remark} As soon as $f$ is not identically zero (i.e. $\nu_f\not={\cal L}_{\T}$), it is possible that $\int_{\T}gd\nu_f\not=0$, while $\int_{\T}g(x)dx=0$. Indeed, there exist an interval $I$ and a real $t$ such that $\nu_f(I)\not=\nu_f(I+t)$, so $g=1_{I}-1_{I+t}$ convenes. In item $ii)$ of the theorem, the condition $g=h-e^{-f}Th$, for some bounded $h$, implies $\int_{\T}gd\nu_{f}=0$. Reciprocally, when $\int_{\T}gd\nu_{f}=0$ and supposing a Diophantine condition on $\theta$ together with a regularity condition on both $f$ and $g$, one can find a bounded $h$ so that $g=h-e^{-f}Th$. For instance, one has the following statement:

\begin{lemme}

$ $

\noindent
Introduce the Diophantine type of $\theta$:

\begin{equation}
\label{type}
\eta(\theta)=\sup\{r\in\R_+^*,~\liminf_{q\rightarrow+\infty} q^r\|q\theta\|=0\}\geq1.
\end{equation}

\noindent
Let $m>\eta(\theta)$ be an integer. Assume that $f\in C^{2m}(\T,\R)$ is centered and $g\in C^m(\T,\R)$ verifies $\int_{\T}g~d\nu_f=0$. Then $g=h-e^{-f}Th$, for some continuous $h$ on $\T$.
\end{lemme}

\noindent
{\it Proof of the lemma:}
 
\noindent
By Arnold \cite{arnold}, cf also Conze-Marco \cite{conzemarco} (Thm 2.1), since $f^{(m)}$ is $C^m$ and $m>\eta(\theta)$, one has $f^{(m)}=v-Tv$, for some continuous $v$. By successive integrations, we have $f=u-Tu$, with $u$ of class $C^m$ and zero mean. Hence $e^{-f}=e^{Tu}/e^u$ and so $\nu_f$ is the measure with density $e^u$ with respect to ${\cal L}_{\T}$. The hypothesis on $g$ is thus $\int (ge^u)(x)dx=0$. As $ge^u$ is of class $C^m$, using one more time \cite{arnold}, we have $ge^u=H-TH$, for a continuous $H$. Finally, $h=e^{-u}H$ is bounded, as continuous on $\T$, and satisfies $g=h-e^{-f}Th$.

\fin

\medskip
In the context of Theorem \ref{thm2} $ii)$, when $\int_{\T}gd\nu_f=0$, for instance when $g=h-e^{-f}Th$ with $h$ bounded (and even simply when $g=0$), transience requires some strongly dissymmetric behaviour between $v_+(n)$ and $w_+(n)$ or between $v_-(n)$ and $w_-(n)$, as $n\rightarrow+\infty$. We build an example in the next proposition. Condition 1) of Theorem \ref{thm2} $ii)$, for example satisfied for $f=1_{[0,1/2)}-1_{[1/2,1)}$ with $x_0=1/4$, prevents this dissymmetry to occur.
\end{remark}

\begin{prop} 
\label{propi}

$ $ 

\noindent
In the context of Theorem \ref{thm2}, there exists $\theta\not\in\Q$ and some BV centered $f$ so that $f=u-Tu$, with $u\geq0$, such that for any bounded $g$, for ${\cal L}_{\T}$-a.-e. $x$, the random walk is transient.
\end{prop}

\subsection{Preliminaries} 

As in \cite{jb2}, we introduce functions $\Phi_{str}(n)$, $\Phi(n)$ and $\Phi_+(n)$ describing the average horizontal macrodispersion of the environment. The last two respectively correspond to $\Phi_u(n)$ and $\Phi_{u,+}(n)$ in \cite{jb2}, Definition 2.3, with $d=1$,  $u=1\in\R_+$ and $\varepsilon_s=m_s$, with the notations of \cite{jb2}.

\begin{defi} 
\label{maindef}

$ $

\noindent
i) The structure function, depending only on the vertical, is defined for $n\geq0$ by:

$$\Phi_{str}(n)=\({n\sum_{-v_-^{-1}(n)\leq k\leq v_+^{-1}(n)}\frac{1}{\rho_k}}\)^{1/2}.$$

\medskip
\noindent
2) For $m,n\geq0$, introduce:

$$\Phi(-m,n)=\({\sum_{-v_-^{-1}(m)\leq k\leq \ell\leq v_+^{-1}(n)}\rho_k\rho_{\ell}\[{\frac{1}{\rho^2_k}+\frac{1}{\rho_{\ell}^2}+\({\sum_{s=k}^{\ell}\frac{\gamma_s\varepsilon_s}{\alpha_s\rho_s}}\)^2}\]}\)^{1/2}.$$

\smallskip
\noindent
For $n\geq 0$, set $\Phi(n)=\Phi(-n,n)$ and $\Phi_{+}(n)=\sqrt{\Phi^2(-n,0)+\Phi^2(0,n)}$.

\end{defi}

\medskip
As in \cite{jb3}, we rectify a misleading point appearing in \cite{jb2}, Definition 2.3 $1)$, where the term ``standard Lebesgue measure'' on the half Euclidean ball $S_+^{d-1}=\{x\in\R^d~|~\|x\|=1,~x_1\geq0\}$ has to be understood as ``uniform probability measure". The following result is extracted from \cite{jb2}, Theorem 2.4, Proposition 2.5 $1)$ and Lemma 6.11.

 \begin{thm}
 \label{theeo}
 
$ $
 
\noindent
i) The random walk is recurrent if and only if $\displaystyle\sum_{n\geq1}\frac{1}{n^2}\frac{(\Phi^{-1}(n))^2}{\Phi_+^{-1}(n)}=+\infty$.

\smallskip
\noindent
ii) The condition $\sum_{n\geq1}1/\Phi(n)<+\infty$ is sufficient for the transience of the random walk. It is equivalent to transience whenever $\Phi\preceq \Phi_+$.

\end{thm}

As a general fact, it is rather directly verified that $\Phi_{str}\preceq\Phi_+\preceq\Phi$. Other general results, fully detailed in \cite{jb3}, section 3.1, are:

\begin{equation}
\label{phiplus}
\Phi_+(n)\asymp\Phi_{str}(n)+\({\sum_{-v_-^{-1}(n)\leq k\leq \ell\leq 0\mbox{ or }0\leq k\leq \ell\leq v_+^{-1}(n)}\rho_k\rho_{\ell}\({\sum_{s=k}^{\ell}\frac{\gamma_s\varepsilon_s}{\alpha_s\rho_s}}\)^2}\)^{1/2},
\end{equation}

\noindent
as well as:

\begin{equation}
\label{phi}\Phi(n)\asymp\Phi_{str}(n)+\({\sum_{-v_-^{-1}(n)\leq k\leq \ell\leq v_+^{-1}(n)}\rho_k\rho_{\ell}\({\sum_{s=k}^{\ell}\frac{\gamma_s\varepsilon_s}{\alpha_s\rho_s}}\)^2}\)^{1/2}.
\end{equation}

\medskip
\noindent
An essential point (recalled in detail in \cite{jb3}, end of Section 3.1) is that the inverse functions $\Phi_{str}^{-1}$, $\Phi_+^{-1}$ and $\Phi^{-1}$ check dominated variation in the sense of Definition \ref{domiv}.

\subsection{Proof of Theorem \ref{thm2} $i)$}

Let us start with a lemma, inspired from \cite{cg} (Proposition 4.2). Introduce the notation :

$$A(n,g,x)=\frac{\sum_{k=0}^ng(T^kx)/\rho_k(x)}{\sum_{k=0}^n1/\rho_k(x)}.$$

\noindent
Recall that $f$ is fixed as in Theorem \ref{thm2}.

\begin{lemme} 
\label{meassu}

$ $

\noindent
i) (Partial reproof of Theorem \ref{cgg}) Let $(x_n)_{n\geq1}$ in $\T$ and $(N_n)$ be an increasing sequence of integers. For $x\in\T$, write $\rho_n(x)=e^{f_n(x)}$, $n\in\Z$. Then any cluster point $\mu$, for the weak-$*$ topology, of the sequence of probability measures on $\T$:

\begin{equation}
\label{mun}
\left(\frac{\sum_{k=0}^{N_n}\delta_{T^kx_n}/\rho_k(x_n)}{\sum_{k=0}^{N_n}1/\rho_k(x_n)}\right)_{n\geq1}
\end{equation}
 
\noindent
is non-atomic and verifies $dT\mu=e^{T^{-1}f}d\mu$. The solution of this last equation is unique.

\medskip
\noindent
ii) Let $g : \T\rightarrow\R$ be BV. Then $A(n,g,x)\rightarrow_{n\rightarrow+\infty}\int_{\T}gd\nu_f$, uniformly in $x$.
\end{lemme}

\noindent
{\it Proof of the lemma:} 

\noindent
As a preliminary remark, for any $(x_n)$ and $(N_n)$, we have, as $n\rightarrow+\infty$:

\begin{equation}
\label{ecras}
\sum_{k=0}^{N_n}1/\rho_k(x_n)\rightarrow+\infty\mbox{ and }(1/\rho_{N_n}(x_n))/(\sum_{k=0}^{N_n}1/\rho_k(x_n))\rightarrow 0.
\end{equation}

\medskip
\noindent
The first point comes from the observation that, independently of $x_n$, $\rho_{q_{\ell}}(x_n)\geq e^{-V(f)}$, for any $\ell\geq1$, as follows from the Denjoy-Koksma inequality \eqref{dk}. Next, the second point can be equivalently rewritten as:

$$\sum_{k=0}^{N_n}e^{-f_k(x_n)+f_{N_n}(x_n)}=\sum_{k=0}^{N_n}e^{T^kf_{N_n-k}(x_n)}=\sum_{k=0}^{N_n}e^{T^{N_n-k}f_{k}(x_n)}\rightarrow+\infty,$$

\medskip
\noindent
for the same reason. 

\medskip
\noindent
$i)$ Call $(\mu_n)_{n\geq1}$ the sequence in \eqref{mun} and consider a cluster point $\mu$ of it for the weak-$*$ topology. For simplicity, we keep the same notations $(x_n)$,$(N_n)$ and assume that $(\mu_n)$ converges to $\mu$. Let $a\in\T$. For any $\delta>0$, we have $\mu((a-\delta,a+\delta))\leq\liminf\mu_n((a-\delta,a+\delta))$. It is therefore enough to show that $\mu_n((a-\delta,a+\delta))$ is arbitrary small for large $n$, for a well-chosen $\delta>0$.

\medskip
Fix an integer $K\geq1$ and take $\delta>0$ so that the intervals $(a-\delta,a+\delta)-k\theta$, $0\leq k\leq 2q_K$, on $\T$ are disjoint. Then:

$$\mu_n((a-\delta,a+\delta))=\frac{\sum_{k=0}^{N_n}1_{(a-\delta,a+\delta)-k\theta}(x_n)/\rho_k(x_n)}{\sum_{k=0}^{N_n}1/\rho_k(x_n)}.$$
 
\medskip
\noindent
Let $0\leq \tau_{1,n}<\cdots<\tau_{L_n,n}\leq N_n$, for some $L_n\geq0$, be the subsequence of $0\leq k\leq N_n$ such that $x_n\in(a-\delta,a+\delta)-k\theta$. If $L_n=0$, we have $\mu_n((a-\delta,a+\delta))=0$. When $L_n\geq 2$, by hypothesis on $\delta$, we have $\tau_{k,n}+q_K<\tau_{k+1,n}$, for $1\leq k<L_n$, and $\tau_{L_n-1,n}+2q_K<\tau_{L_n,n}$. Using the Denjoy-Koksma inequality \eqref{dk}, giving $1/\rho_{k\pm q_l}(x_n)\geq e^{-V(f)}/\rho_{k}(x_n)$, we obtain, when $L_n\geq2$:

\begin{eqnarray}
\sum_{k=0}^{N_n}\frac{1}{\rho_k(x_n)}&\geq&\sum_{1\leq k\leq L_n-1}\sum_{0\leq \ell\leq K}\frac{1}{\rho_{\tau_{k,n}+q_{\ell}}(x_n)}+\sum_{0\leq \ell\leq K}\frac{1}{\rho_{\tau_{L_n,n}-q_{\ell}}(x_n)}\nonumber\\
&\geq&\sum_{1\leq k\leq L_n}(1/\rho_{\tau_{k,n}}(x_n))\sum_{0\leq \ell\leq K}e^{-V(f)}=(K+1)e^{-V(f)}\sum_{1\leq k\leq L_n}(1/\rho_{\tau_{k,n}}(x_n)).\nonumber\end{eqnarray}

\smallskip
\noindent
When $L_n=1$, noticing that either $\tau_{L_n,n}>q_K$ or $\tau_{L_n,n}+q_K< N_n$, as soon as $n$ is large enough, the same reasoning holds and we obtain the same equality. Hence, always:

$$\mu_n((a-\delta,a+\delta))\leq\frac{\sum_{1\leq k\leq L_n}1/\rho_{\tau_{k,n}}(x_n)}{(K+1)e^{-V(f)}\sum_{1\leq k\leq L_n}(1/\rho_{\tau_{k,n}}(x_n))}=\frac{e^{V(f)}}{K+1}.$$
 
\medskip
\noindent
This can be made arbitrary small, when choosing $K$ large enough. Hence $\mu$ is non-atomic. 

\medskip
Next, for any continuous $h :\T\rightarrow\R$, $A(N_{n},h,x_{n})\rightarrow \int_{\T}hd\mu$. Since $\mu$ is non-atomic, this holds for any $h$ continuous except at countably many points and in particular if $h$ is BV. Since $f$ is BV, $e^{-f}$ is also BV. Thus for any continuous $h$, $A(N_{n},e^{-f}Th,x_{n})\rightarrow \int_{\T}e^{-f}Thd\mu$. It now follows from \eqref{ecras}, that for any continuous $h$:

$$\int_{\T}e^{-f}Thd\mu=\int_{\T}hd\mu,$$

\medskip
\noindent
giving $dT\mu=e^{T^{-1}f}d\mu$. Thus $\mu$ solves the equation $dT\mu=e^{T^{-1}f}d\mu$. For unicity of the solution, see \cite{cg} (Theorem 5.6). This solution is hence non-atomic. This completes the proof of this point.

\medskip
\noindent
$ii)$ If the result is not true, using that $|A(n,g,x)|\leq\|g\|_{\infty}$, there exists $a\in\R$, $(x_n)$ in $\T$ and $N_n\rightarrow+\infty$, such that $A(N_n,g,x_n)\rightarrow a\not=\int_{\T}gd\nu_f$. By compacity of the weak-$*$ topology, for some subsequence $(N_{\zeta(n)})$ of $(N_n)$ and $(x_{\zeta(n)})$ of $(x_n)$, the sequence of measures:

$$\frac{\sum_{k=0}^{N_{\zeta(n)}}\delta_{T^kx_{\zeta(n)}}/\rho_k(x_{\zeta(n)})}{\sum_{k=0}^{N_{\zeta(n)}}1/\rho_k(x_{\zeta(n)})},n\geq1,$$

\smallskip
\noindent
converges to some probability $\mu$ on $\T$. By $i)$, $dT\mu=e^{T^{-1}f}d\mu$ and so $\mu=\nu_f$, by unicity. As $\nu_f$ is non-atomic, $A(N_{\zeta(n)},g,x_{\zeta(n)})\rightarrow \int_{\T}gd\nu_f$, contradicting $A(N_n,g,x_n)\rightarrow a\not=\int_{\T}gd\nu_f$. This ends the proof of the lemma. 

\fin

\bigskip
We turn to the proof of Theorem \ref{thm2} $i)$, fixing BV functions $f$ and $g$, with $f$ centered and $\int_{\T}gd\nu_f\not=0$. Let $x\in\T$ and $\rho_n=\rho_n(x)=e^{f_n(x)}$, $n\in\Z$, as before. Observe first that (cf \cite{jb3}, section 3.1, 2), for the first inequality):

$$\Phi_{str}(n)\preceq \({\sum_{-v_-^{-1}(n)\leq k\leq \ell\leq v_+^{-1}(n)}\left(\frac{\rho_k}{\rho_{\ell}}+\frac{\rho_{\ell}}{\rho_k}\right)}\)^{1/2}\preceq\({\sum_{-v_-^{-1}(n)\leq k\leq \ell\leq v_+^{-1}(n)}\rho_k\rho_{\ell}\({\sum_{s=k}^{\ell}\frac{1}{\rho_s}}\)^2}\)^{1/2},$$

\medskip
\noindent
considering, for the second inequality, only in the last inside sum the terms for $s=k$ and $s=\ell$. Introduce now the following function $\Psi$, essentially corresponding to $\Phi$ when $g=1$ (``essentially", because, by \eqref{phi} and the previous inequality, the definition can be simplified when $g=1$):

\begin{equation}
\label{psii}
{\Psi}(n)=\({\sum_{-v_-^{-1}(n)\leq k\leq \ell\leq v_+^{-1}(n)}\rho_k\rho_{\ell}\({\sum_{s=k}^{\ell}\frac{1}{\rho_s}}\)^2}\)^{1/2}.
\end{equation}

\medskip
\noindent
Notice now that (cf again \cite{jb3}, section 3.1, 2), for the first line), using at the end that $g$ is bounded:

\begin{eqnarray}
\Phi(n)&\asymp&\Phi_{str}(n)+\({\sum_{-v_-^{-1}(n)\leq k\leq \ell\leq v_+^{-1}(n)}\rho_k\rho_{\ell}\({\sum_{s=k}^{\ell}\frac{\gamma_s\varepsilon_s}{\alpha_s\rho_s}}\)^2}\)^{1/2}\nonumber\\
&=&\Phi_{str}(n)+\({\sum_{-v_-^{-1}(n)\leq k\leq \ell\leq v_+^{-1}(n)}\rho_k\rho_{\ell}\({\sum_{s=k}^{\ell}\frac{g(T^sx)}{\rho_s}}\)^2}\)^{1/2}\preceq{\Psi}(n).\nonumber
\end{eqnarray}

\medskip
\noindent
We now prove the reverse inequality. Using Lemma \ref{meassu} $ii)$, let first $M\geq1$ be such that for $n\geq M$ and all $x\in\T$: 

$$\left|\frac{\sum_{k=0}^ng(T^kx)/\rho_k(x)}{\sum_{k=0}^n1/\rho_k(x)}-\int_{\T}gd\nu_f\right| \leq \left|\int_{\T}gd\nu_f\right|/2.$$

\medskip
\noindent
This gives the inequality $|\sum_{s=k}^{\ell}g(T^sx)/\rho_s(x)|\geq(1/2)|\int_{\T}gd\nu_f|\times \sum_{s=k}^{\ell}1/\rho_s(x)$, whenever $\ell-k>M$. Consequently:

\begin{eqnarray}
\label{notzero}
\Phi(n)&\geq&\({\sum_{-v_-^{-1}(n)\leq k\leq \ell\leq v_+^{-1}(n),\ell-k>M}\rho_k\rho_{\ell}\({\sum_{s=k}^{\ell}\frac{g(T^sx)}{\rho_s}}\)^2}\)^{1/2}\nonumber\\
&\geq&\frac{1}{2}\left|\int_{\T}gd\nu_f\right|\({\sum_{-v_-^{-1}(n)\leq k\leq \ell\leq v_+^{-1}(n),\ell-k>M}\rho_k\rho_{\ell}\({\sum_{s=k}^{\ell}\frac{1}{\rho_s}}\)^2}\)^{1/2}.\end{eqnarray}

\medskip
\noindent
Observe now that:

$$\({\sum_{-v_-^{-1}(n)\leq k\leq \ell\leq v_+^{-1}(n),\ell-k\leq M}\rho_k\rho_{\ell}\({\sum_{s=k}^{\ell}\frac{g(T^sx)}{\rho_s}}\)^2}\)^{1/2}\leq \|g\|_{\infty}\({\sum_{-v_-^{-1}(n)\leq k\leq \ell\leq v_+^{-1}(n),\ell-k\leq M}\rho_k\rho_{\ell}\({\sum_{s=k}^{\ell}\frac{1}{\rho_s}}\)^2}\)^{1/2}.$$

\medskip
\noindent
As $M$ is fixed and $\eta\leq\rho_k(y)/\rho_{k+1}(y)\leq 1/\eta$, for any $k\in\Z$ and $y\in\T$, where $\eta$ comes from Hypothesis \ref{hypo}, we get when $0\leq \ell-k\leq M$:

$$\rho_k\rho_{\ell}\({\sum_{s=k}^{\ell}\frac{1}{\rho_s}}\)^2\asymp 1.$$ 

\smallskip
\noindent
Therefore, for some constant $C>0$, depending on $M$:

\begin{eqnarray}
\label{racine}
\({\sum_{-v_-^{-1}(n)\leq k\leq \ell\leq v_+^{-1}(n),\ell-k\leq M}\rho_k\rho_{\ell}\({\sum_{s=k}^{\ell}\frac{1}{\rho_s}}\)^2}\)^{1/2}&\leq& C\({\sum_{-v_-^{-1}(n)\leq k\leq v_+^{-1}(n)}1}\)^{1/2}\nonumber\\
&\preceq&\sqrt{v_-^{-1}(n)+v_+^{-1}(n)}.
\end{eqnarray}

\medskip
\noindent
We now show that $v_+^{-1}(n)/\Psi^2(n)\rightarrow0$. Indeed, by \eqref{psii}:

\begin{eqnarray}
\Psi^2(n)\geq\sum_{0\leq k\leq \ell\leq v_+^{-1}(n)}\rho_k\rho_{\ell}\({\sum_{s=k}^{\ell}\frac{1}{\rho_s}}\)^2&\succeq&\sum_{0\leq k\leq \ell\leq v_+^{-1}(n)}\({\frac{\rho_k}{\rho_{\ell}}+\frac{\rho_{\ell}}{\rho_k}}\)\nonumber\\
&\succeq&\sum_{0\leq k\leq v_+^{-1}(n)}\rho_k\sum_{0\leq k\leq v_+^{-1}(n)}1/\rho_{k}\succeq(v_+^{-1}(n))^2,\nonumber
\end{eqnarray}

\smallskip
\noindent
by the Cauchy-Schwarz inequality in the final step. Thus $\Psi^2(n)/v_+^{-1}(n)\succeq v_+^{-1}(n)\rightarrow+\infty$, as $n\rightarrow+\infty$. In the same way, $\Psi^2(n)/v_-^{-1}(n)\rightarrow+\infty$, as $n\rightarrow+\infty$. By \eqref{racine}:

$$\({\sum_{-v_-^{-1}(n)\leq k\leq \ell\leq v_+^{-1}(n),\ell-k\leq M}\rho_k\rho_{\ell}\({\sum_{s=k}^{\ell}\frac{1}{\rho_s}}\)^2}\)^{1/2}=o(\Psi(n)).$$

\medskip
\noindent
From \eqref{notzero} and \eqref{racine}, we deduce that $\Psi\preceq\Phi$ and thus finally $\Psi\asymp\Phi$. The same argumentation shows that $\Phi_+\asymp\Psi_+$, where $\Psi$ ``essentially" corresponds to $g=1$ and is defined by: 

\begin{equation}
\label{psii+}
\Psi_+(n)=\({\sum_{-v_-^{-1}(n)\leq k\leq \ell\leq v_+^{-1}(n),k\ell\geq0}\rho_k\rho_{\ell}\({\sum_{s=k}^{\ell}\frac{1}{\rho_s}}\)^2}\)^{1/2}.
\end{equation}

\medskip
By Theorem \ref{theeo} $i)$, we are therefore left to proving:

\begin{equation}
\label{final}
\sum_{n\geq1}\frac{1}{n^2}\frac{(\Psi^{-1}(n))^2}{\Psi_+^{-1}(n)}<+\infty.
\end{equation}

\medskip
\noindent
We shall give two proofs of \eqref{final}. The first one just a reinterpretation. Consider another random walk, this time defined by (changing only $g$ and keeping the same BV function $f$ and $x\in\T$):

$$\mu_n=\delta_{+1},~\varepsilon_n=1,~\alpha_n=\gamma_n,~\beta_n/\alpha_n=e^{f(T^{n-1}x)},~n\in\Z.$$

\smallskip
\noindent
Since $\gamma_n\varepsilon_n/\alpha_n=1$, this case corresponds to $g=1$. As previously indicated (before \eqref{psii}), the functions $\tilde{\Phi}$ and $\tilde{\Phi}_+$ in this case verify $\tilde{\Phi}\asymp\Psi$ and $\tilde{\Phi}_+\asymp\Psi_+$. By Theorem \ref{theeo} $i)$, condition \eqref{final} is just the transience criterion of this new random walk. The latter being obviously transient (a.-s., the horizontal coordinate tends monotonically to $+\infty$), the condition in \eqref{final} is verified and this ends the first proof.

\medskip
One may be interested in showing directly the convergence of the series in \eqref{final}. We now furnish the argument. This may help in the future to manipulate the recurrence criterion for studying other examples. Introduce the following functions:

\begin{equation}
\label{psi++}
\Psi_{++}(n)=\({\sum_{0\leq k\leq \ell\leq v_+^{-1}(n)}\rho_k\rho_{\ell}\({\sum_{s=k}^{\ell}\frac{1}{\rho_s}}\)^2}\)^{1/2},
\end{equation}

\smallskip
\noindent 
as well as:

$$\Psi_{+-}(n)=\({\sum_{-v_-^{-1}(n)\leq k\leq \ell\leq 0}\rho_k\rho_{\ell}\({\sum_{s=k}^{\ell}\frac{1}{\rho_s}}\)^2}\)^{1/2}.$$

\medskip
\noindent
We shall next use repeatedly properties like $\max(a,b)\asymp a+b$, $\sqrt{a+b}\asymp\sqrt{a}+\sqrt{b}$, etc, for $a,b\geq0$. From the definition of $\Psi_+$ given in \eqref{psii+}, we have:

$$\Psi_+(n)\asymp \Psi_{++}(n)+\Psi_{+-}(n).$$

\medskip
\noindent
Hence, $\Psi^{-1}_+(n)\asymp \min\{\Psi^{-1}_{++}(n),\Psi^{-1}_{+-}(n)\}$, thus furnishing:

$$\frac{1}{\Psi^{-1}_+(n)}\asymp \frac{1}{\Psi^{-1}_{++}(n)}+\frac{1}{\Psi^{-1}_{+-}(n)}.$$

\medskip
\noindent
In order to prove \eqref{final}, we thus have to show the two convergences:

\begin{equation}
\label{final1}
\sum_{n\geq1}\frac{1}{n^2}\frac{(\Psi^{-1}(n))^2}{\Psi_{++}^{-1}(n)}<+\infty\mbox{ and }\sum_{n\geq1}\frac{1}{n^2}\frac{(\Psi^{-1}(n))^2}{\Psi_{+-}^{-1}(n)}<+\infty.
\end{equation}

\medskip
\noindent
We establish the first one, the case of the second one being similar. 

\medskip
First, using Hypothesis \ref{hypo}, for some $c>1$ depending only on $\eta$, we have $\sum_{0\leq k\leq n+1}\rho_k\leq c\sum_{0\leq k\leq n}\rho_k$, for all $n\geq0$. Hence, for large $n>0$, uniformly on $x\in\T$, using \eqref{sand}:

\begin{equation}
\label{minori}
n\geq\sum_{0\leq \ell\leq v_+^{-1}(n)}\rho_{\ell}\geq \sum_{0\leq \ell\leq v_+^{-1}(n)+1}\rho_{\ell}/c\geq n/c.
\end{equation}

\smallskip
\noindent
Thus, for large $n>0$ (uniformly on $x\in\T$):

\begin{equation}
\label{entredeux}
\sum_{v_+^{-1}(n/c^2)\leq \ell\leq v_+^{-1}(n)}\rho_{\ell}\asymp n,
\end{equation}

\medskip
\noindent
We now have, using \eqref{psii}:

\begin{equation}
\label{pz}
\Psi(n)\geq\({\sum_{-v_-^{-1}(n)\leq k\leq 0\leq \ell\leq v_+^{-1}(n)}\rho_k\rho_{\ell}\({\sum_{s=0}^{\ell}\frac{1}{\rho_s}}\)^2}\)^{1/2}.
\end{equation}

\smallskip
\noindent
The variables $k$ and $\ell$ are now independent on the right hand side. Let us define:

$$\zeta(n)=\sqrt{n}\({\sum_{0\leq \ell\leq v_+^{-1}(n)}\rho_{\ell}\({\sum_{s=0}^{\ell}\frac{1}{\rho_s}}\)^2}\)^{1/2}.$$

\smallskip
\noindent
Obviously, by \eqref{minori}:

\begin{equation}
\label{des}
\zeta(n)\preceq n\sum_{s=0}^{v_+^{-1}(n)}\frac{1}{\rho_s}.
\end{equation}

\smallskip
\noindent
Analogously to \eqref{minori}, we have $\sum_{-v_-^{-1}(n)\leq k\leq 0}\rho_k\asymp n$, so we get, by \eqref{pz} and \eqref{entredeux}:

\begin{eqnarray}
\label{tr}
\Psi(n)\succeq\zeta(n)&\geq&\sqrt{n}\({\sum_{v_+^{-1}(n/c^2)\leq \ell\leq v_+^{-1}(n)}\rho_{\ell}\({\sum_{s=0}^{\ell}\frac{1}{\rho_s}}\)^2}\)^{1/2}\nonumber\\
&\succeq& n\sum_{s=0}^{v_+^{-1}(n/c^2)}\frac{1}{\rho_s}\geq (n/c^2)\sum_{s=0}^{v_+^{-1}(n/c^2)}\frac{1}{\rho_s}.\end{eqnarray}

\smallskip
\noindent
Set $F(n)=\sum_{k=0}^{v_+^{-1}(n)}1/\rho_k$ and $G(n)=nF(n)$. The first inequality in \eqref{tr} gives $\Psi^{-1}(n)\preceq\zeta^{-1}(n)$. Moreover, the last inequalities in \eqref{tr} also provide, together with \eqref{des}:

$$\zeta^{-1}(n)\asymp G^{-1}(n).$$

\smallskip
\noindent
In order to establish the first part of \eqref{final1}, it is thus sufficient to show that:

\begin{equation}
\label{final2}
\sum_{n\geq1}\(\frac{G^{-1}(n)}{n}\)^2\frac{1}{\Psi^{-1}_{++}(n)}<+\infty.
\end{equation}

\medskip
For the analysis of $\Psi_{++}$, we fix an integer $K>c^2$, where the constant $c>1$ appears in \eqref{minori}. Define now for any integer $u\geq0$ the quantity:

\begin{equation}
\label{au}
A_u=\sum_{v_+^{-1}(K^u)<k\leq v_+^{-1}(K^{u+1})}1/\rho_k.
\end{equation}

\smallskip
\noindent
Starting from the definition \eqref{psi++} of $\Psi_{++}$:

\begin{eqnarray}
\Psi_{++}(K^n)&=&\(\sum_{0\leq k\leq \ell\leq v_+^{-1}(K^n)}\rho_k\rho_l\(\sum_{k\leq s\leq\ell}1/\rho_s\)^2\)^{1/2}\nonumber\\
&\preceq&\(\sum_{v_+^{-1}(K^0)< k\leq \ell\leq v_+^{-1}(K^n)}\rho_k\rho_l\(\sum_{k\leq s\leq\ell}1/\rho_s\)^2\)^{1/2}.\nonumber
\end{eqnarray}

\smallskip
\noindent
Proceeding as for \eqref{entredeux}, we also have $\sum_{v_+^{-1}(K^u)<k\leq v_+^{-1}(K^{u+1})}\rho_k\asymp K^u$. Using this, we continue:

\begin{eqnarray}
\Psi_{++}(K^n)&\preceq&\(\sum_{0\leq u\leq v\leq n-1}\sum_{\underset{v_+^{-1}(K^v)<\ell\leq v_+^{-1}(K^{v+1}),k\leq\ell}{v_+^{-1}(K^u)<k\leq v_+^{-1}(K^{u+1})}}\rho_k\rho_l\(\sum_{v_+^{-1}(K^u)<s\leq v_+^{-1}(K^{v+1})}1/\rho_s\)^2\)^{1/2}\nonumber\\
&\preceq&\(\sum_{0\leq u\leq v\leq n-1}K^uK^v(A_u+\cdots+A_v)^2\)^{1/2}.\nonumber\end{eqnarray}

\medskip
\noindent
We arrive at:

\begin{eqnarray}
\label{zab1}
\Psi_{++}(K^n)&\preceq& \sum_{0\leq u\leq v< n}K^{u/2}K^{v/2}(A_u+\cdots+A_v)\nonumber\\
&\preceq& \sum_{0\leq l< n}A_l\sum_{0\leq u\leq \ell}K^{u/2}\sum_{\ell\leq v<n}K^{v/2}\preceq K^{n/2}\sum_{0\leq \ell< n}A_{\ell}K^{\ell/2}.\end{eqnarray}

\medskip
\noindent
Let $N\geq1$ be an integer. We have, by \eqref{zab1}:

\begin{eqnarray}
\label{zab3}
\sum_{n=0}^N\frac{\Psi_{++}(K^n)}{K^n}\preceq\sum_{0\leq n\leq N}\frac{1}{K^{n/2}}\sum_{0\leq \ell<n}A_{\ell}K^{\ell/2}=\sum_{0\leq \ell< N}A_{\ell}\sum_{\ell<n\leq N}\frac{K^{\ell/2}}{K^{n/2}}\asymp\sum_{0\leq \ell< N}A_{\ell}\asymp F(K^{N}).
\end{eqnarray}

\medskip
\noindent
Now, \eqref{zab3} furnishes, still for $N\geq1$:

\begin{eqnarray}
\label{fk}
\sum_{0\leq k\leq\Psi_{++}(K^N)}\frac{1}{\Psi^{-1}_{++}(k)}&\preceq&\sum_{0\leq n<N}\sum_{\Psi_{++}(K^n)<k\leq\Psi_{++}(K^{n+1})}\frac{1}{\Psi^{-1}_{++}(k)}\nonumber\\
&\preceq&\sum_{0\leq n<N}\frac{\Psi_{++}(K^{n+1})}{K^n}\preceq F(K^N).
\end{eqnarray}

\medskip
\noindent
Let us define: 

$$Z(n)=\sum_{0\leq k\leq n}1/\Psi^{-1}_{++}(k).$$

\smallskip
\noindent
We extend the notation to real $x>0$, by $Z(x)=Z(\lfloor x\rfloor)$, using the floor function. Idem for $\Psi_{++}$. The last inequality \eqref{fk} thus says that for $n\geq1$, taking $N$ so that $K^{N-1}\leq n<K^N$:

$$Z(\Psi_{++}(n/K))\preceq F(n)=\frac{G(n)}{n}.$$

\medskip
\noindent
In particular, using \eqref{sand}:

$$Z(\Psi_{++}(G^{-1}(n)/K))\preceq \frac{G(G^{-1}(n))}{G^{-1}(n)}\preceq \frac{n}{G^{-1}(n)}.$$

\medskip
\noindent
Next, if $n\leq \Psi_{++}(G^{-1}(n)/K)$, then $Z(n)\preceq n/G^{-1}(n)$. Otherwise, the last inequality gives:

\begin{eqnarray}
Z(n)&\leq& Z(\Psi_{++}(G^{-1}(n)/K))+\sum_{\Psi_{++}(G^{-1}(n)/K)<k\leq n}1/\Psi^{-1}_{++}(k)\nonumber\\
&\leq& \frac{n}{G^{-1}(n)}+\frac{n}{\Psi^{-1}_{++}(\Psi_{++}(G^{-1}(n)/K))}\preceq\frac{n}{G^{-1}(n)}.\nonumber\end{eqnarray}

\medskip
\noindent 
Finally, using the last inequality and the definition of $Z(n)$, we show \eqref{final2}:

\begin{eqnarray}
\sum_{n\geq1}\(\frac{G^{-1}(n)}{n}\)^2\frac{1}{\Psi^{-1}_{++}(n)}\preceq\sum_{n\geq1}\(\frac{1}{Z(n)}\)^2\frac{1}{\Psi^{-1}_{++}(n)}&\leq&\sum_{n\geq1}\frac{1}{Z(n-1)Z(n)}\frac{1}{\Psi^{-1}_{++}(n)}\nonumber\\
&\leq&\sum_{n\geq1}\(\frac{1}{Z(n-1)}-\frac{1}{Z(n)}\)<+\infty.\nonumber
\end{eqnarray}

\medskip
\noindent
This concludes the second proof of \eqref{final} and of Theorem \ref{thm2} $i)$.

\fin

\medskip
\noindent
\begin{remark}
The previous proof in fact shows that, in complete generality, the condition $\Phi_{str}^2\preceq\Phi$ implies the transience of the random walk. 
\end{remark}

\subsection{Proof of Theorem \ref{thm2} $ii)$}

In this section, $f$ and $g$ are BV functions, with $f$ centered. We suppose that $g=h-e^{-f}Th$, for a bounded $h$. As $T^sg(x)=\gamma_s\varepsilon_s/\alpha_s$ and $T^s(e^{-f}Th)/\rho_s=(T^{s+1}h)/\rho_{s+1}$, we first have:

\begin{eqnarray}
\({\sum_{-v_-^{-1}(n)\leq k\leq \ell\leq v_+^{-1}(n)}\rho_k\rho_{\ell}\({\sum_{s=k}^{\ell}\frac{\gamma_s\varepsilon_s}{\alpha_s\rho_s}}\)^2}\)^{1/2}&=&\({\sum_{-v_-^{-1}(n)\leq k\leq \ell\leq v_+^{-1}(n)}\rho_k\rho_{\ell}\({\sum_{s=k}^{\ell}\frac{T^sg}{\rho_s}}\)^2}\)^{1/2}\nonumber\\
&\preceq&\({\sum_{-v_-^{-1}(n)\leq k\leq \ell\leq v_+^{-1}(n)}\rho_k\rho_{\ell}\({\sum_{s=k}^{\ell}\({\frac{T^sh}{\rho_s}-\frac{T^{s+1}h}{\rho_{s+1}}}\)}\)^2}\)^{1/2}\nonumber\\
&\preceq&\({\sum_{-v_-^{-1}(n)\leq k\leq \ell\leq v_+^{-1}(n)}\rho_k\rho_{\ell}\({1/\rho_k^2+1\rho^2_{l}}\)}\)^{1/2}\nonumber\\
&\asymp& \Phi_{str}(n).\nonumber
\end{eqnarray}

\medskip
\noindent
See for example \cite{jb3}, section 3.1, 2), for details of the last step. Next, we deduce by \eqref{phi} that $\Phi(n)\asymp \Phi_+(n)\asymp\Phi_{str}(n)$. Using Theorem \ref{theeo} $ii)$, the recurrence of the random walk is now equivalent to the divergence of $\sum_{n\geq1}1/\Phi_{str}(n)$, or, using \eqref{entredeux} and the analogous version for $v_-$, of (cf Definition \ref{4.1}):

$$\sum_{n\geq1}\frac{1}{\sqrt{n(w_+\circ v_+^{-1}(n)+w_-\circ v_-^{-1}(n))}}.$$

\medskip
\noindent
Because of the monotonicity of the general term of the previous series, by usual condensation, this is equivalent, for any fixed $K>1$, to showing the divergence of:

\begin{equation}
\label{divergence}
\sum_{n\geq1}\frac{\sqrt{K^n}}{\sqrt{w_+\circ v_+^{-1}(K^n)+w_-\circ v_-^{-1}(K^n)}}.
\end{equation}

\medskip
\noindent
1) Suppose first that $f=u-Tu$, with $e^{u}\in L^1({\cal L}_{\T})$. Then, by the Law of Large Numbers (for the second step), a.-e., as $n\rightarrow+\infty$: 

\begin{equation}
\label{lln}
w_+(n)(x)=\sum_{k=0}^ne^{-u(x)+T^ku(x)}\sim ne^{-u(x)}\int_{\T}e^{u(y)}dy.
\end{equation}

\smallskip
\noindent
In the same way, a.-e., $w_-(n)(x)$ is linear. For $v_+$, using again the Law of Large Numbers, for a positive, but maybe non-integrable, function, there is a.-e. some $\kappa(x)>0$ such that:

$$v_+(n)(x)=e^{u(x)}\sum_{k=0}^ne^{-T^ku(x)}\geq n\kappa(x)\mbox{, for }n\geq1.$$

\smallskip
\noindent
The same property is true, a.-e., for $v_-(n)$, as $n\rightarrow+\infty$. Hence, a.-e., there is some $c(x)>0$ so that $v_+^{-1}(n)\leq c(x)n$ and $v_-^{-1}(n)\leq c(x)n$, for large $n$. We obtain, a.-e., for large $n>0$:

$$w_+\circ v_+^{-1}(n)\leq w_+(c(x)n)\preceq n.$$

\smallskip
\noindent
Idem, $w_-\circ v_-^{-1}(n)\preceq n$. These make the general term in \eqref{divergence} not go to zero, so the series diverges.

\medskip
\noindent
2) Suppose that (instead of the $L^1$-condition) for some $x_0\in\T$, then $f(x+x_0)=f(x_0-x)$, for a.-e. $x\in\T$. Using the denominators $(q_n)$ of the convergents of the angle $\theta$, one has:

\begin{eqnarray}
\label{tutu}
\sum_{0\leq k\leq q_n}e^{f_k(x)}=\sum_{0\leq k\leq q_n}e^{f_{q_n-k}(x)}&=&\sum_{0\leq k\leq q_n}e^{f_{-k}(x)+T^{-k}f_{q_n}(x)}\asymp\sum_{0\leq k\leq q_n}e^{f_{-k}(x)},
\end{eqnarray}

\medskip
\noindent
using the Denjoy-Koksma's inequality \eqref{dk}. As a result, one obtains that $v_+(q_n)\asymp v_-(q_n)$, as $n\rightarrow+\infty$, uniformly on $x\in\T$. For the same reason:

\begin{equation}
\label{ww+}
w_+(q_n)\asymp w_-(q_n),
\end{equation}

\smallskip
\noindent
as $n\rightarrow+\infty$, also uniformly on $x\in\T$. 

\medskip
\noindent
From $v_+(q_n)\asymp v_-(q_n)$, independently on $x\in\T$, we can now fix some large $K>1$ and $p_0$ such that for any $n\geq1$, there exists $p$ with $K^p\leq v_+(q_n)\leq K^{p+p_0}$ and $K^p\leq v_-(q_n)\leq K^{p+p_0}$. This gives $v_+^{-1}(K^p),v_-^{-1}(K^p)\leq q_n$. In \eqref{divergence}, the term corresponding to $p$ verifies (uniformly on $x\in\T$):

\begin{eqnarray}
\label{v+w+}
\frac{\sqrt{K^p}}{\sqrt{w_+\circ v_+^{-1}(K^p)+w_-\circ v_-^{-1}(K^p)}}\geq\frac{\sqrt{K^{-p_0}v_+(q_n)}}{\sqrt{w_+(q_n)+w_-(q_n)}}\asymp\frac{\sqrt{v_+(q_n)}}{\sqrt{w_+(q_n)}},
\end{eqnarray}

\smallskip
\noindent
using \eqref{ww+} for the last step. Next, immediately from the definition of the model, the set $\{x\in\T\mbox{, the random walk is transient}\}$ is measurable and $T$-invariant, hence has Lebesgue measure zero or one, by ergodicity of $(\T,T,{\cal L}_{\T})$. If the random walk were transient for a.-e. $x$, then, by \eqref{v+w+} and the convergence of the series in \eqref{divergence}, for a.-e. $x$:

\begin{equation}
\label{tz}(v_+(q_n)/w_+(q_n))(x_0+x)\rightarrow0\mbox{ and }(v_+(q_n)/w_+(q_n))(x_0-x)\rightarrow0,
\end{equation}

\medskip
\noindent
as $n\rightarrow+\infty$. However, using \eqref{ww+} and the symmetry assumption in the final step, we can write:

\begin{eqnarray}
\frac{v_+(q_n)}{w_+(q_n)}(x_0+x)\asymp\frac{v_+(q_n)}{w_-(q_n)}(x_0+x) =\frac{\sum_{0\leq k\leq q_n}e^{f_k(x_0+x)}}{\sum_{0\leq k\leq q_n}e^{-f_{-k}(x_0+x)}}\asymp\frac{\sum_{0\leq k\leq q_n}e^{f_k(x_0+x)}}{\sum_{0\leq k\leq q_n}e^{f_{k}(x_0-x)}}.\nonumber
\end{eqnarray}

\medskip
\noindent
As a result, for a.-e. $x\in\T$:

$$\frac{v_+(q_n)}{w_+(q_n)}(x_0+x)\asymp \frac{w_+(q_n)}{v_+(q_n)}(x_0-x).$$

\smallskip
\noindent
By \eqref{tz}, the left hand side goes to $0$, as $n\rightarrow+\infty$, whereas the right hand side goes to $+\infty$. This contradiction completes the proof of Theorem \ref{thm2} $ii)$.

\fin

\medskip
\noindent
\begin{remark}In a similar way, but without the symmetry assumption, suppose that $f(x)=u(x)-u(x+y)$, for some BV function $u$ and some parameter $y\in\T$. Let us show that for a.-e. $(x,y)\in\T^2$ the random walk is recurrent. Indeed, in the previous proof, part 2), if ever transience holds for some $(x,y)$, then, by \eqref{v+w+} and the convergence of the series in \eqref{divergence}:

$$v_+(n)/w_+(n)=\(\sum_{k=0}^{q_n}e^{u_k(x)-u_k(x+y)}\)/\(\sum_{k=0}^{q_n}e^{-u_k(x)+u_k(x+y)}\)\rightarrow0.$$

\medskip
\noindent
Now, the set of $(x,y)\in\T^2$ verifying this property is clearly invariant by the joint action of $T\times Id$ and $Id\times T$ on $\T^2$, which is ergodic. Hence this set has measure 0 or 1. If this is 1, one obtains that for a.-e. $(x,y)$, as $n\rightarrow+\infty$:

$$\(\sum_{k=0}^{q_n}e^{u_k(x)-u_k(y)}\)/\(\sum_{k=0}^{q_n}e^{-u_k(x)+u_k(y)}\)\rightarrow0.$$

\medskip
\noindent
This is impossible again, when reversing the roles of $x$ and $y$. We thus have recurrence for a.-e. $(x,y)\in\T^2$. Rather generally, in Theorem \ref{thm2} $ii)$, it would be interesting if the symmetry assumption 2) could be dropped. This raises the question, for $f:\T\rightarrow\R$, BV and centered, of understanding the a.-e. behaviour, as $n\rightarrow\infty$, of ratios of the form $\sum_{k=0}^{n}e^{f_k(x)}/\sum_{k=0}^{n}e^{-f_k(x)}$.

\medskip
\noindent
Also in Theorem \ref{thm2} $ii)$, supposing only $\int_{\T} gd\nu_f=0$ for $g$ (in place of $g=h-e^{-f}Th$, with $h$ bounded) requires to find sharp upper-bounds on sums of the form:

$$\sum_{k=0}^ne^{-f_k(x)}T^kg(x).$$

\end{remark}

\subsection{Proof of Proposition \ref{propi}}

By Theorem \ref{theeo} $ii)$, to prove transience for the random walk, it is enough to show the convergence of $\sum_{n\geq1}1/\Phi_{str}(n)$. By definition of $\Phi_{str}$, it is sufficient to establish that:

\begin{equation}
\label{ridi}
\sum_{n\geq1}1/\sqrt{nw_+\circ v^{-1}_+(n)}<+\infty.
\end{equation}

\smallskip
\noindent
We choose $f$ in the form $f=u-Tu$, with $u\geq 0$. Then $e^{-u}$ is integrable, so $v_+(n)$ is a.-e. linear, as $n\rightarrow+\infty$, by the Law of Large Numbers, as for \eqref{lln}. As a consequence, we obtain that for any $x\in\T$, $v_+^{-1}(n)\asymp n$, as $n\rightarrow+\infty$. 

\medskip
\noindent
Setting $U=e^u$, we have $w_+(n)\sim e^{-u(x)}U_n(x)$, so in order to obtain \eqref{ridi}, it is enough to show that for a.-e. $x$:

\begin{equation}
\label{thefinal}
\sum_{n\geq1}1/\sqrt{nU_n}<+\infty.
\end{equation}

\smallskip
Let us build $u$ (and $f=u-Tu$). Let the rotation angle $\theta\not\in\Q$ be defined by the partial quotients $a_m=m^6$, $m\geq1$. Introduce $h_{B,\Delta}(x)=B(1-\vert x\vert/\Delta)_+$, for $\Delta>0$, $B>0$. It is a piecewise linear pick function of height $B$ and width $\Delta$, centered at zero. 

\medskip
\noindent
Let $(q_m)$ be the denominators of the convergents of $\theta$. For $m\geq1$, set $h^m=h_{B_m,\Delta_m}$, with:

\begin{equation}
\label{Del}
\Delta_m=1/(m^2q_m),~~~~~B_m=m^2/q_m.
\end{equation}

\medskip
\noindent
We first define $f=\sum_{m\geq1}f^m$, where:

$$f^m=\sum_{k=0}^{q_m-1}T^{-k}(h^m-T^{q_m}h^m).$$

\medskip
\noindent
For large $m$, the $(k\theta)_{0\leq k<q_m}$ are approximately equally spaced and the sum in the definition of $f^m$ involves functions with disjoint supports. For $m\geq1$, $f^m$ is centered and, using \eqref{Del} and \eqref{kh1}:

$$V(f^m)\leq q_mV(h^m-T^{q_m}h^m)\leq C q_m(B_m/\Delta_m)\Vert q_m\theta\Vert\leq C/m^2.$$

\medskip
\noindent
As a result, $f$ is BV and centered. We now check that $f$ is a.-e. equal to some $u-Tu$. For $m\geq1$, $f^m=u^m-Tu^m$, with:

\begin{equation}
\label{um}
u^m=\sum_{k=0}^{q_m-1}T^{-k}\sum_{\ell=0}^{q_m-1}T^{\ell}h^m=\sum_{\vert \ell\vert< q_m}(q_m-\vert \ell\vert)T^{\ell}h^m.
\end{equation}

\medskip
\noindent
The Lebesgue measure of the support of $u^m$ is $\leq 2q_m\Delta_m$. As $\sum_{m\geq1}q_m\Delta_m<+\infty$, by the first lemma of Borel-Cantelli, a.-e. $x\in\T$ belongs to the support of $u^m$, for only finitely many $m$. Hence $u=\sum_{m\geq1}u^m$ is well-defined a.-e. and $f=u-Tu$, a.-e..

\medskip
The Diophantine type $\eta(\theta)$ of $\theta$, defined \eqref{type}, also equal $\limsup \log q_{n+1}/\log q_n$. Hence, here $\eta(\theta)=1$, as $a_m=m^6$, $m\geq1$. For $x\in \T$ and $r>0$, let $\tau_r(x)=\min\{n\geq 1,\Vert T^nx\Vert<r\}$. By a result of Kim and Marmi \cite{kimarmi}, for a.-e. $x$:

\begin{equation}
\label{kim}
\lim_{r\rightarrow0}\frac{\log \tau_r(x)}{-\log r}=1.
\end{equation}

\medskip
\noindent
Recall that $U=e^u$ and let $s_m=\tau_{\Delta_m/2}(x)$, $m\geq1$. For a.-e. $x\in\T$, decompose now:

\begin{equation}
\label{kimi}
\sum_{n> s_1}\frac{1}{\sqrt{nU_n}}=\sum_{m\geq1}\sum_{s_m< n\leq s_{m+1}}\frac{1}{\sqrt{nU_n}}\preceq\sum_{m\geq1}\sqrt{s_{m+1}}/\sqrt{U_{s_m+1}}.
\end{equation}

\medskip
\noindent
We next have, using \eqref{um} and $h^m(T^{s_m}x)\geq B_m/2$:

\begin{equation}
\label{den1}
U_{s_m+1}=e^{u_{s_m+1}(x)}\geq e^{u(T^{s_m}x)}\geq e^{u^m(T^{s_m}x)}\geq e^{q_mh^m(T^{s_m}x)}\geq e^{m^2/2}.
\end{equation}

\medskip
\noindent
Also, a.-e., using \eqref{kim}, for large $m$, $s_m\leq(2/\Delta_m)^2\preceq (m^2q_m)^2$. Next, the denominators of the convergents of $\theta$ verify $q_m=m^6q_{m-1}+q_{m-2}\leq (2m)^6q_{m-1}$, so brutally $q_m\leq (2m)^{6m}\leq e^{cm\log m}$. Hence $s_{m+1}=O(e^{3cm\log m})$. Together with \eqref{den1}, we deduce that the series in \eqref{kimi} is finite. This ends the proof of the proposition. 

\fin

\bigskip
As a final remark, rather generally, when $f$ is not an additive coboundary, then the asymptotic behaviour of $(\log(v_+(n)/w_+(n)))_{n\geq1}$ is somehow that of a classical random walk in $\R$. 

\begin{lemme}

$ $ 

\noindent
Let $f\in {\cal L}^1(\T,\R)$, centered, and not a.-e. equal to $u-Tu$, for some measurable $u$. Define $v_+(n)=v_+(n)(x)$ and $w_+(n)=w_+(n)(x)$, $n\geq1$, as in Definition \ref{4.1}. Then one of the following three situations occur:

\medskip
\noindent
1) $v_+(n)/w_+(n)\rightarrow+\infty$, a.s., as $n\rightarrow+\infty$.

\noindent
2) $v_+(n)/w_+(n)\rightarrow0$, a.s., as $n\rightarrow+\infty$.

\noindent
3) $\limsup v_+(n)/w_+(n)\rightarrow+\infty$, a.s., and $\liminf v_+(n)/w_+(n)\rightarrow0$, a.s., as $n\rightarrow+\infty$.

\end{lemme}

\medskip
\noindent
{\it Proof of the lemma :} 

\noindent
As a first point, $v_+(n)\rightarrow+\infty$, a.-e., since, a.-e., the random walk $(f_k)_{k\geq0}$ is recurrent, as $k\rightarrow+\infty$, since $f$ is integrable and centered. We have $v_+(n)(x)\sim e^{f(x)}v_+(n-1)(Tx)$ and $w_+(n)(x)\sim e^{-f(x)}w_+(n-1)(Tx)$, so the following set is $T$-invariant:

$$\{x\in\T,~\limsup v_+(n)(x)/w_+(n)(x)<+\infty\}.$$ 

\smallskip
\noindent
It hence has Lebesgue measure 0 or 1, by ergodicity. If this measure is 1, we can a.-e. define $\psi(x)=\limsup_{n\rightarrow+\infty}(v_+(n)(x)/w_+(n)(x))$. The opening remark on equivalents yields:

$$\psi(x)=e^{2f(x)}\psi(Tx).$$

\smallskip
\noindent
Now, the set $\{\psi(x)>0\}$ is $T$-invariant and thus again has measure 0 or 1. If this measure is 1, one has $f=(\log\psi)/2-(\log T\psi)/2$, contrary to the hypothesis. Hence the set has measure 0. 

\medskip
Finally, $\limsup v_+(n)/w_+(n)=+\infty$, a.-s., or $v_+(n)/w_+(n)\rightarrow0$, a.-s., as $n\rightarrow+\infty$. Symmetrically, $\limsup w_+(n)/v_+(n)=+\infty$, a.-s., or $w_+(n)/v_+(n)\rightarrow0$, a.-s.. Intersecting the possibilities, we obtain the three cases given in the statement of the lemma.

\fin

\bigskip
\noindent
{\bf{Acknowledgments.}} We thank Yves Derriennic for inviting to revisit these questions and the referee for detailed comments.

\providecommand{\bysame}{\leavevmode\hbox to3em{\hrulefill}\thinspace}

\bigskip
{\small{\sc{Univ Paris Est Creteil, CNRS, LAMA, F-94010 Creteil, France\\
Univ Gustave Eiffel, LAMA, F-77447 Marne-la-Vall\'ee, France}}}

\it{E-mail address~:} {\sf julien.bremont@u-pec.fr}

\end{document}